\documentclass{amsart}

\usepackage{graphicx}

\usepackage{amsmath, amssymb, array, latexsym, hyperref}

\usepackage{fullpage}

\newtheorem{theorem}{Theorem}[section]
\newtheorem{lemma}[theorem]{Lemma}
\newtheorem{proposition}[theorem]{Proposition}
\newtheorem{corollary}[theorem]{Corollary}

\theoremstyle{definition}

\newtheorem{examples}[theorem]{Examples}

\theoremstyle{remark}
\newtheorem{remark}[theorem]{Remark}
\newtheorem{remarks}[theorem]{Remarks}
\newtheorem{note}[theorem]{Note}

\numberwithin{equation}{section}

\newcommand{\prend}{$\hfill \Box$}
\newcommand{\ls}{\leqslant}
\newcommand{\gr}{\geqslant}

\begin{document}

\title{Variance estimates and almost Euclidean structure}

\author{Grigoris Paouris}
\address{Department of Mathematics, Mailstop 3368, Texas A\&M University, College Station, TX, 77843-3368.}  
\email{grigorios.paouris@gmail.com}
\thanks{The first author was supported by the NSF CAREER-1151711 grant.} 

\author{Petros Valettas}
\address{Mathematics Department, University of Missouri, Columbia, MO, 65211.}
\email{valettasp@missouri.edu}
\thanks{The second author was supported by the NSF grant DMS-1612936.}



\begin{abstract}
We introduce and initiate the study of new parameters associated with any norm and any log-concave measure 
on $\mathbb R^n$, which provide sharp distributional inequalities. In the Gaussian context this investigation sheds light to 
the importance of the statistical measures of dispersion of the norm in connection with the local structure of the ambient space. 
As a byproduct of our study, we provide a short proof of Dvoretzky's theorem which not only supports the aforementioned significance 
but also complements the classical probabilistic formulation.
\end{abstract}

\maketitle



\section{Introduction}

The main focus of this note is to establish new distributional inequalities for convex functions with respect to log-concave 
measures. The new ingredient in these inequalities is that the controlling parameter is the variance rather than the Lipschitz constant or
 some moment of the length of the gradient of the function. The study of these inequalities has been motivated by the need to quantify 
 efficiently the almost Euclidean structure of high-dimensional normed spaces. In particular, we would like to understand
 the dependence on $\varepsilon$ in the almost isometric version of Dvoretzky's theorem \cite{Dvo}. In fact in this note, we are interested in 
 the randomized Dvoretzky's theorem (as has been established in the seminal work of Milman \cite{Mil1}) which states that for any $\varepsilon\in (0,1)$ 
 there exists a constant $c(\varepsilon)>0$ with the following property: for any normed space 
 $X=(\mathbb R^n, \|\cdot\|)$ the random (with respect to the Haar measure on the Grassmannian) $k$-dimensional subspace $F$ 
 satisfies
 \begin{align} \label{eq:local-osc}
(1-\varepsilon) M \|z\|_2 \ls\|z\| \ls (1+\varepsilon) M \|z\|_2
 \end{align} for all $z\in F$ as long as $k\ls c(\varepsilon)k(X)$, where $M$ is the average of the norm on the unit Euclidean sphere with respect to the uniform probability measure (see \S2 for the related definitions). The parameter $k(X)$ is referred to as the {\it critical dimension} (or 
the {\it Dvoretzky number}) of $X$ and is given by $k(X)=n(M(X)/b(X))^2$, where $b(X)=\max\{\|\theta\| : \|\theta\|_2=1\}$. We will write 
$k(X,\varepsilon)$ for the maximal
dimension for which \eqref{eq:local-osc} holds with probability at least $1/2$. Thus, Milman's formulation yields $k(X,\varepsilon)\gr c(\varepsilon)k(X)$. (For historical remarks, required background material and further extensions the reader may consult \cite{MS}, \cite{Pis_book}, \cite{TJ} and \cite{AGM}). Concerning the dependence on $\varepsilon$, let us mention that Milman's proof
yields $c(\varepsilon)\simeq \varepsilon^2/\log\frac{1}{\varepsilon}$. This was subsequently improved by Gordon in \cite{Go} showing
that one can always have $c(\varepsilon) \simeq \varepsilon^2$ and an alternative approach was presented by Schechtman 
in \cite{Sch1}. Hence, one has $k(X,\varepsilon)\gr c\varepsilon^2 k(X)$. The latter asymptotic formula is optimal up to universal constants as the example of the $\ell_1^n$ shows. However, there are (several) examples of spaces which show that the function $c(\varepsilon)$ can be significantly improved.

The investigation of the randomized Dvoretzky's theorem on the almost isometric level for special cases of normed 
spaces was initiated by the works of Schechtman \cite{Sch-cube} and Tikhomirov \cite{Tik-cube} who determined asymptotically the dimension $k(\ell_\infty^n, \varepsilon)$. Their estimate $k(\ell_\infty^n,\varepsilon)\simeq \varepsilon |\log \varepsilon|^{-1} \log n$ is better than 
$\varepsilon^2 \log n$ as was predicted by the optimal form of Milman's formula due Gordon's \cite{Go}. 
This study has been extended by Zinn and the authors in \cite{PVZ} for all $\ell_p^n$ spaces. 
Once again the dimension $k(\ell_p^n,\varepsilon)$ is much larger
than the estimate derived by the classical proof when $2<p<\infty$. The proof is based on the fact that the $\ell_p$ norm is 
much more concentrated on Gauss' space than the usual ``concentration of measure" suggests. In that case we will say 
that the norm is ``overconcentrated".

This phenomenon is not only observed in the $\ell_p^n$ spaces, but is apparent in many other cases.
In \cite{PV-dvo-L_p} we have shown that for all $2<p<\infty$ and for any $n$-dimensional subspace $X$ of $L_p$,
in Lewis' position \cite{Lew}, either $X$ has almost spherical sections of proportional dimension (that is $k(X)\simeq n$) or $X$ is 
overconcentrated. 
Furthermore, among all $n$-dimensional subspaces $X$ of $L_p, \; 2<p<\infty$ the worst $k(X,\varepsilon)$ occurs for the $\ell_p^n$'s.
Recently Tikhomirov \cite{Tik-1-unc} proved that for every 1-unconditional
$n$-dimensional normed space $X$, in $\ell$-position (see e.g. \cite{TJ} for the definition of the $\ell$-position), either $k(X)$ is polynomial with respect to $n$ or $X$
is overconcentrated. His approach shows that the worst $k(X,\varepsilon)$ in this class of spaces is attained for the $\ell_\infty^n$.

In some of the above cases, the new observed phenomenon is due to superconcentration (as defined by Chatterjee in \cite{Cha}) and 
can be quantified by employing Talagrand's $L_1-L_2$ bound for the gaussian measure (see e.g. \cite{Tal-russo} and \cite{CorLed}). 
That was the crucial tool in the investigation in
\cite{PVZ} as well as in \cite{Tik-1-unc}. However, the superconcentration is not the only reason that causes the norm of a space to be more concentrated than expected. All the aforementioned cases share the common feature 
that ${\rm Var}\|Z\| \ll {\rm Lip}(\|\cdot\|)^2$ where $Z$ is
an $n$-dimensional standard gaussian vector (recall that the classical gaussian concentration yields 
the bound ${\rm Var}\|Z\| \ls {\rm Lip}(\|\cdot\|)^2$).

Our aim in this note is to put on display the importance of the (normalized) variance in the study of the almost Euclidean structure in 
high-dimensional normed spaces and to initiate a systematic investigation of the concentration properties of convex functions with
respect to this measure of dispersion. This investigation has two main directions. 
Firstly, we show that this new parameter can be used to provide a short proof of Dvoretzky's theorem which can be 
viewed as the probabilistic and quantitative version of the topological proof 
due to Figiel \cite{F} and Szankowski's analytic proof from \cite{Sz}. Further study of this parameter 
is also considered and is compared with the
classical Dvoretzky number. 

Secondly, we study distributional inequalities in the context of log-concave measures. We show that (one-sided) 
deviation estimates, where the variance is involved, can be proved for any pair of a convex function and a log-concave 
probability measure on Euclidean space, by extending a machinery developed in \cite{PVsd}.

The rest of the paper is organized as follows: In Section 2 we fix the notation and we recall standard background material. 
In Section 3 we present a simple proof that every infinite dimensional Banach space contains $\ell_2^n$'s uniformly. The local version of our approach uses
the parameter of the normalized variance and yields optimal dependence on the size of the distortion. In Section 4 we study further the
aforementioned parameter and we compare it with the classical Dvoretzky number. Our approach uses 
concentration estimates for the mean width of random projections of any convex body in $\mathbb R^n$. We postpone the somewhat more 
systematic study of this topic to \S6. In Section 5 we study,
in the general context of log-concave measures small deviation and small ball estimates for norms whose tightness is quantified in terms
of the parameter of the normalized variance. Finally, in Section 6 we study independently
distributional inequalities for the mean width of random projections of any convex body in $\mathbb R^n$ in terms of the Haar measure on
the Grassmannian.

\medskip
\section{Notation and background material}

We work in $\mathbb R^n$ equipped with a Euclidean structure $\langle \cdot, \cdot \rangle$. The Euclidean norm is given by
$\|x\|_2 = \sqrt{ \langle x,x \rangle}$, $x\in \mathbb R^n$. For the $(n-1)$-dimensional Euclidean sphere we write 
$S^{n-1}= \{ x\in \mathbb R^n : \langle x, x \rangle =1 \}$. Let also $\sigma$ be the uniform probability measure on $S^{n-1}$.

The orthogonal group on $\mathbb R^n$ is denoted by $O(n)$ and consists of all matrices which preserve the angles, i.e. 
$\langle Ux, Uy \rangle =\langle x, y\rangle$ for all $x,y\in \mathbb R^n$. Thus,
\begin{align*}
O(n)=\{U \in \mathbb R^{n\times n} : UU^\ast= I\}.
\end{align*} The action of $O(n)$ on itself generates the Haar probability measure which we will denote by $\mu_n$.

Let $SO^\pm(n)$ be the collection of all matrices which are volume preserving,
\begin{align*}
SO^\pm(n)= \{U\in O(n) : \det U= \pm 1 \}.
\end{align*} Similarly, we may define the Haar probability measure on the special orthogonal group $SO(n)$. 
This is nothing more than the restriction of $\mu_n$ in $SO(n)$ (see e.g. \cite{Me} for background details), thus we will still denote 
it by $\mu_n$. Note that although $SO^-(n)$ is not a group itself we will refer to the restricted measure, with some abuse 
of terminology, as the Haar measure on $SO^-(n)$ since it is still invariant (within $SO^-(n)$) under the action of $SO^+(n)$. 
One may check that if $U$ is Haar distributed in $SO^+(n)$ and $M$ is a fixed matrix in $SO^-(n)$, then $UM$ is Haar distributed
in $SO^-(n)$. The Hilbert-Schmidt norm of a matrix $T$ is denoted by $\|T\|_{\rm HS}$ and the operator norm of $T$ with $\|T\|_{\rm op}$.

The Grassmann space $G_{n,k}$ is the set of all $k$-dimensional subspaces of $\mathbb R^n$. We consider the Haar measure 
$\nu_{n,k}$ on $G_{n,k}$ which is invariant under the action of the orthogonal group $O(n)$. 
For an arbitrary norm on $\mathbb R^n$ we write $\|\cdot \|$ and for the normed spaces $(\mathbb R^n, \|\cdot\|)$ we will 
use the letters $X,Y$.


The random variables will be denoted by $\zeta, \eta, \xi, \ldots$. For random vectors on $\mathbb R^n$, distributed according to some law $\mu$, we write $Z,W$ or $W=(w_1,\ldots,w_n)$. 
For any random variable $\xi$ on some probability space $(\Omega, \mathcal A, \mathbb P)$ with 
$\xi\in L_2(\Omega)$ and $\mathbb E \xi\neq 0$ we define the parameter $\beta(\xi)$ as follows:
\begin{align}
\beta(\xi) := \frac{ {\rm Var (\xi)}}{(\mathbb E \xi)^2},
\end{align} where $\mathbb E $ stands for the expectation and ${\rm Var}$ for the variance. 
For any normed space $X=(\mathbb R^n , \|\cdot\|)$ and for any Borel probability measure $\mu$ on $\mathbb R^n$ 
we define the parameter:
\begin{align} \label{def-beta}
\beta_\mu(X) := \beta(\|Z\|),
\end{align} where $Z$ is a random vector distributed according to $\mu$. If the prescribed measure is the 
Gaussian we will often omit the subscript. 

Recall the well known fact that if $\xi$ is a random variable in some probability space $(\Omega, \mathcal A, \mathbb P)$ with 
$\mathbb E\xi^2<\infty$ and $m={\rm med}(\xi)$ is a median of $\xi$, then
\begin{align*}
| \mathbb E \xi -m | \leqslant \mathbb E |\xi -m| \leqslant \sqrt{ {\rm Var}( \xi ) }.
\end{align*} If $\mathbb E \xi \neq 0$ then we readily get:
\begin{align*}
\left| \frac{m}{\mathbb E\xi} -1\right| \ls \sqrt{ \beta(\xi) }.
\end{align*} This says that the parameter $\beta$ quantifies how close are the measures of central tendency, median and expectation.

Let us note that for the probability space $(\mathbb R^n,\mu)$, where $\mu$ is a log-concave probability measure and the variable $\xi$ 
is $\|Z\|$, where $\|\cdot\|$ is a norm on $\mathbb R^n$ and $Z$ is distributed according to $\mu$, Borell's lemma \cite{Bor} 
(see also \S 5) implies $m\simeq \mathbb E\|Z\|$. Thus, if we are only interested in the order of magnitude of the parameter $\beta$, 
the expectation can be sufficiently replaced by the median. For general facts concerned with gaussian measures we refer the reader 
to \cite{Bog}.


A convex body $K$ in $\mathbb R^n$ is a compact, convex set with non-empty interior. The convex bodies will be denoted with
$A, K, L, \ldots$. A convex body $K$ is said to be symmetric when $x\in K$ if and only if $-x\in K$. For any symmetric convex
body $K$ we write $\|\cdot\|_K$ for the norm (gauge function) induced by $K$. Therefore, very frequently we identify notationally the normed
space, whose norm is generated by a convex body, with the convex body itself.
The volume (i.e. the Lebesgue measure) of a convex body 
$K$ in $\mathbb R^n$ is denoted by $|K|$. We define the circumradius of $K$, as the number 
$R(K) := \max_{x\in B_2^n}\|x\|_K$. That is the smallest $R>0$ for which $K \subseteq R B_2^n$. For any symmetric convex 
body $K$ in $\mathbb R^n$ we define the following averages:
\begin{align*}
J_q (K):= \left(\frac{1}{|K|} \int_K \|x\|_2^q \, dx\right)^{1/q}, \quad -n<q\neq 0.
\end{align*}

The next lemma is essentially from \cite{Kl1}. The idea of the proof goes
back to Rudelson (see \cite{Kl-sub}):

\begin{lemma}\label{lem: stability of moms}
Let $K$ be a symmetric convex body in $\mathbb R^n$. Then, for all $q>0$ we have:
\begin{align*}
J_q (K) \gr a_{n,q}^{-1} R(K), \quad a_{n,q}^{-q}:=\frac{q}{2}B(q,n+1).
\end{align*} 
\end{lemma} 

\medskip

\begin{remarks} (i) Let us note that by using standard asymptotic estimates for the Beta function we have:
\begin{align}
a_{n,q} \ls \exp\left( \frac{cn}{q} \log \left( \frac{eq}{n} \right)\right)
\end{align} for all $q\gr n$.

\smallskip

\noindent (ii) The symmetry is not essential for the proof. 
\end{remarks}

\noindent Next, for any symmetric convex body $K$ in $\mathbb R^n$ we define:
\begin{align*}
{\rm vrad}(K) := \left( \frac{|K|}{|B_2^n|} \right)^{1/n}, \quad M_q^q(K):= \int_{S^{n-1}} \|\theta \|_K^q \, d\sigma(\theta), \quad q\neq 0.
\end{align*}

The proof of the next lemma can be found in \cite{Kl1}. 

\begin{lemma} \label{lem:identity}
Let $K$ be a symmetric convex body in $\mathbb R^n$. Then, for any $p> -n$ we have:
\begin{align*}
{\rm vrad}(K)^n J_p^p(K) = \frac{n}{n+p} M_{-(n+p)}^{-(n+p)}(K).
\end{align*}
\end{lemma}

We denote by $b(K)=\max_{\theta\in S^{n-1}}\|\theta\|_K$. Note that $r(K)=1/b(K)$ is the maximal radius of the centered inscribed ball in $K$. 
The polar body $K^\circ$ of a convex body $K$ is defined as $K^\circ= \{y : \langle x, y\rangle \ls 1 \; \forall x\in K\}$. For any 
symmetric convex body $K$ the support function of $K$ at $\theta\in S^{n-1}$ is the half width of the body in direction $\theta$, i.e.
$h_K(\theta) =\max\{ \langle x,\theta \rangle : x\in K\}$. Note that $h_K$ is the dual norm of $\|\cdot\|_K$. One may check that $b(K)=R(K^\circ)$ and $R(K)r(K^\circ)=1$. The Banach-Mazur distance between two isomorphic normed spaces $X,Y$ is denoted by $d(X,Y)$ with
\begin{align*}
d(X,Y):= \inf\{ \|T\|\cdot \|T^{-1}\| : T:X\to Y, \; {\rm linear \, isomorphism} \}.
\end{align*}
The reader may consult the monographs \cite{MS,Pis_book,TJ,AGM} for a detailed background information on the local theory of normed spaces.

Throughout the text we make use of $C,c,C_1,c_1,\ldots$ for positive absolute constants whose values may change
from line to line. We also introduce the notation $Q_1 \lesssim Q_2$ for any two quantities $Q_1$ and $Q_2$ 
(which may depend on dimension or some geometric parameter of the normed space or the convex body) if there exists absolute constant $C>0$ such that $Q_1 \ls CQ_2$. We write
$Q_1\simeq Q_2$ is $Q_1\lesssim Q_2$ and $Q_2\lesssim Q_1$. 
We write the above signs with a subscript $\lesssim_p, \simeq_p$, if the involved constant depends on the parameter $p$.

\medskip
\section{A probabilisitic proof of Dvoretzky's theorem without concentration: The variance method}

In this section we provide a probabilistic proof that $\ell_2^n$'s embed uniformly into any infinite dimensional Banach 
space $X$ without utilizing the concentration of measure. What is crucial in our approach is the use of the new critical 
parameter $\beta (X)$ defined in \eqref{def-beta}. We investigate this parameter further in \S 3.

In order to prove the result we will use some standard lemmas such as the Dvoretzky-Rogers lemma and a net argument on the 
sphere which are rather standard in all probabilistic proofs. However, even if our approach is more elementary, it
leads to better dependence on $\varepsilon$ in several interesting cases (see \S 3.3). Our proposed local form of Dvoretzky's theorem 
reads as follows.

\begin{theorem} \label{thm:w-dvo}
For every $\varepsilon \in (0,1)$ and for any finite-dimensional normed space $X$, there exists 
$k$-dimensional subspace $F$ of $X$ with $k\gr c\log \frac{1}{\beta} / \log \frac{1}{\varepsilon}$, such that $d(F, \ell_2^k)<1+\varepsilon$,
where $\beta=\beta (X)$.
\end{theorem}

Later this formulation will be used to compare the new parameter $\beta (X)$ with the Dvoretzky number $k(X)$ of the normed space $X=(\mathbb R^n, \|\cdot\|)$. For the proof we will need the following standard lemma 
(whose proof is included for the sake of completeness).

\begin{lemma} \label{lem:approx-net} Let $X$ be a $k$-dimensional normed space.

\smallskip

\noindent {\rm (i)} For any $\delta\in (0,1)$ there exists a $\delta$-net $\mathcal N$ 
on $S_X$ with ${\rm card}(\mathcal N) \ls \left (1+ \frac{2}{\delta} \right)^k$.

\smallskip

\noindent {\rm (ii)} Let $Y$ be a normed space and let $T: X\to Y$ be a linear mapping with the property:
\begin{align*}
1-\varepsilon \ls \|Tz\| \ls 1+\varepsilon
\end{align*} for all $z$ in a $\delta$-net $\mathcal N$ of $S_X$ where $0<\delta, \varepsilon<1$. Then, we have:
\begin{align*}
\frac{1-\varepsilon -2\delta}{1-\delta}  \ls \|T\theta\| \ls \frac{1+\varepsilon}{1-\delta} ,
\end{align*} for all $\theta \in S_X$. 
\end{lemma}

\noindent {\it Proof.} The proof of the first assertion can be found in \cite{MS}. Let us proceed with (ii). Fix $\theta\in S_X$. Then, there 
exists $z\in \mathcal N$ with $\|z-\theta\|<\delta$. Thus, we may write:
\begin{align} \label{eq: 2.3}
\|T\theta\| \ls \delta\|T\| +\|Tz\| \ls \delta \|T\| +(1+\varepsilon).
\end{align} Since $\theta$ was arbitrary, it follows that $ \|T\| \ls \delta \|T\| +(1+\varepsilon)$, or equivalently
\begin{align} \label{eq: 2.4}
\|T\| \ls (1-\delta)^{-1}(1+\varepsilon) .
\end{align} Plug this back in \eqref{eq: 2.3}
we obtain the right-hand side estimate. For the left hand side we argue as follows:
\begin{align*}
\|T\theta\| \gr \|Tz\|-\|T(\theta-z)\| \gr (1-\varepsilon)  -\delta \|T\| \gr (1-\varepsilon)- \frac{\delta (1+\varepsilon) }{1-\delta} ,
\end{align*} by the estimate \eqref{eq: 2.4}. \prend

\bigskip

Now we are ready to prove the aforementioned form of Dvoretzky's theorem:

\medskip

\noindent {\it Proof of Theorem \ref{thm:w-dvo}.} Let $\{g_{ij}(\omega)\}_{i,j=1}^{n,k}$ be i.i.d. standard Gaussian random 
variables in some probability space $(\Omega, \mathcal A, \mathbb P)$. We consider the random Gaussian operator: $G:\mathbb \ell_2^k \to X$ with 
\begin{align*}
G_\omega z= \sum_{i=1}^n \sum_{j=1}^k g_{ij}(\omega) z_j e_i , \quad z=(z_1,\ldots,z_k )\in \ell_2^k, \; \omega\in \Omega.
\end{align*} Note that for fixed $\theta\in S^{k-1}$ if we apply Chebyshev's inequality we obtain:
\begin{align*}
\mathbb P \left( \big| \|G\theta\| -\mathbb E\|Z\| \big| > \varepsilon \mathbb E\|Z\| \right) = 
\mathbb P \left( \big| \|Z\| -\mathbb E\|Z\| \big| > \varepsilon \mathbb E\|Z\| \right) \ls \beta / \varepsilon^2,
\end{align*} for $\varepsilon>0$ and $Z$ standard Gaussian random vector. 
Now fix $\varepsilon\in (0,1/3)$ and employ Lemma \ref{lem:approx-net} with $\delta= \varepsilon /2$, to get:
\begin{align*}
\mathbb P \left( \exists z\in \mathcal N : \big| \|Gz\|- \mathbb E\|Z\| \big|> \varepsilon \mathbb E\|Z\| \right) \ls \left( \frac{6}{\varepsilon}\right)^k \frac{\beta}{\varepsilon^2} 
\ls (6/ \varepsilon)^{3k} \beta.
\end{align*} Therefore, as long as $(6 / \varepsilon)^{3k} \beta<1/2$ or $k\ls c \log \frac{1}{\beta}/ \log \frac{1}{\varepsilon}$, 
there exists $\omega \in \Omega$ with the property:
\begin{align*}
\big| \|G_\omega z\|- \mathbb E\|Z\| \big| \ls \varepsilon \mathbb E\|Z\|,
\end{align*} for all $z\in \mathcal N$. Fix the operator $G=G_\omega$ and let $T:\ell_2^k \to X$ with $T: =(\mathbb E\|Z\|)^{-1} G$. 
Then, by Lemma \ref{lem:approx-net} and the choice of $\delta$ we conclude that:
\begin{align*}
1-2\varepsilon< \frac{1-\varepsilon -2\delta}{1-\delta}  \ls \|T\theta\| \ls \frac{1+\varepsilon}{1-\delta} <1 +2\varepsilon ,
\end{align*} for all $\theta \in S^{n-1}$. The result follows; for $F:=T(\ell_2^k)$ we get $d(F,\ell_2^k) <1+11\varepsilon$. \prend

\subsection{Quantitative form of Theorem \ref{thm:w-dvo}.} 

Using John's position \cite{Jo} and the classical Dvoretzky-Rogers 
Lemma from \cite{DR} we can show the following estimate:

\begin{proposition} \label{prop: beta-est} Let $X=(\mathbb R^n, \| \cdot\|)$ be a finite-dimensional normed space 
and assume that the Euclidean ball $B_2^n$ is the ellipsoid of maximal volume inside $B_X$. Then, we have:
\begin{align} \label{eq: beta-est}
\beta (X) \ls \frac{C}{\log n}.
\end{align}
\end{proposition}

Therefore, taking into account Proposition \ref{prop: beta-est} we readily get the finite representability of $\ell_2$ in 
any infinite dimensional Banach space $X$:

\begin{corollary}
For every $\varepsilon \in (0,1)$ and for any $n$-dimensional normed space $X$ there exists $k\gr c \log\log n/\log \frac{1}{\varepsilon}$ and $k$-dimensional subspace $F$ of $X$ with $d(F, \ell_2^k)<1+\varepsilon$. 
\end{corollary}

\noindent {\it Proof of Proposition \ref{prop: beta-est}. (Sketch).} In order to prove the above estimate we need the next lemma from \cite{DR}:

\begin{lemma} \label{lem:dr-lem}
Let $X=(\mathbb R^n, \|\cdot\|)$ be $n$-dimensional normed space and let $B_2^n$ be the ellipsoid of maximal volume inside $B_X$.
Then, there exist orthonormal vectors $v_1,\ldots, v_m$ with $m \simeq n$ such that $\|v_j\|\gr 1/2$ for $j=1,2,\ldots,m$.
\end{lemma}

A proof of this result can be found in \cite{MS}. Using Lemma \ref{lem:dr-lem} and an averaging argument over signs (see \cite{MS} for details)
we arrive at the following estimate:

\smallskip

\noindent {\it Claim 1.} If $B_2^n$ is the ellipsoid of maximal volume inside $B_X$ we have:
\begin{align*}
\mathbb E\|Z\| \gr c\sqrt{\log n}, \quad Z\sim N({\bf 0},I_n).
\end{align*}

One more ingredient is needed:

\smallskip

\noindent {\it Claim 2.} If $\|\cdot\|$ is a norm on $\mathbb R^n$ and $b= \max_{\theta\in S^{n-1}} \|\theta\|$, then
\begin{align*}
{\rm Var} \|Z\| \ls b^2.
\end{align*} 

\noindent For the proof of Claim 2 we may employ the (Gaussian) Poincar\'e inequality:
\begin{align} \label{eq:Poin}
{\rm Var}[f(Z)] \ls \mathbb E \big \|  \nabla f(Z) \big\|_2^2.
\end{align} The fact that $\langle \nabla \|x\|,v \rangle \ls b$ for all $x\in \mathbb R^n$ and $v\in S^{n-1}$ proves the claim.

The assertion follows if we combine the previous two claims. \prend

\subsection{Probabilistic estimate} 

The argument in Theorem \ref{thm:w-dvo} in fact implies that for
$k\ls c \log (1/\beta) / \log (1/\varepsilon)$ the {\it random} $k$-dimensional subspace $F$ of $X$ is $(1+\varepsilon)$-Euclidean 
with probability $>1-e^{-3k}$. Let us first recall the definition: 

Let $X=(\mathbb R^n, \|\cdot\|)$ be a normed space. A subspace $F$ of $X$ is said to be
{\it $(1+\varepsilon)$-spherical} if 
\begin{align*} 
\max_{z\in S_F}\|z\| / \min_{z\in S_F} \|z\|< 1+\varepsilon.
\end{align*}

In order to verify the aforementioned probabilistic statement, note that the proof of Theorem \ref{thm:w-dvo} yields the following estimate:
\begin{align*}
P \left( \big| \|G\theta\| -\mathbb E\|Z\| \big| \ls \varepsilon \mathbb E\|Z\|, \, \forall \theta\in S^{k-1} \right) \gr 1-e^{-3k},
\end{align*} provided that $k\ls \frac{1}{3} \log (1/\beta) / \log (6e/\varepsilon)$. Note that:
\begin{align*}
\left \{ \big| \|G\theta\| -\mathbb E\|Z\| \big| \ls \varepsilon \mathbb E\|Z\|, \; \forall \theta \in S^{k-1} \right\} 
\subseteq \left \{ \max_\theta \|G\theta\| / \min_\theta \|G\theta\|  < \frac{1+\varepsilon}{1-\varepsilon}  \right \}\equiv \mathcal A_1
\end{align*} and similarly 
\begin{align*}
\left \{ \forall \theta \in S^{k-1} \, ; \,\big| \|G\theta\|_2 -\mathbb E\|Z\|_2 \big| \ls \varepsilon \mathbb E\|Z\|_2\right\} 
\subseteq \left \{ \max_\theta \|G\theta\|_2 / \min_\theta \|G\theta\|_2  < \frac{1+\varepsilon}{1-\varepsilon}  \right \} \equiv \mathcal A_2.
\end{align*} Furthermore, the event $\mathcal A_1 \cap \mathcal A_2$ satisfies:
\begin{align*}
\mathcal A_1 \cap \mathcal A_2 \subseteq \left\{ G(\mathbb R^k) \, {\rm is} \, \left(\frac{1+\varepsilon}{1-\varepsilon}\right)^2-{\rm spherical} \right\}.
\end{align*} Finally, recall that $P( \omega : G_\omega(\mathbb R^k) \in \mathcal B) = \nu_{n,k}(F \in G_{n,k} : F\in \mathcal B)$ for any Borel set 
$\mathcal B$ in $G_{n,k}$ and that 
\begin{align*}
P\left(\left \{ \big| \|G\theta\|_2 -\mathbb E\|Z\|_2 \big| \ls \varepsilon \mathbb E\|Z\|_2, \; \forall \theta\in S^{k-1} \right\} \right) >1-Ce^{-c\varepsilon^2 n}.
\end{align*} Thus, we get:
\begin{align*}
\nu_{n,k} \left( \left\{ F \, {\rm is} \,  \left( \frac{1+\varepsilon}{1-\varepsilon} \right)^2-{\rm spherical} \right\} \right) \gr P(\mathcal A_1 \cap \mathcal A_2)
&\gr 1-e^{-3k}-Ce^{-c\varepsilon^2 n} \\
&\gr 1-C'e^{-c' k},
\end{align*} since $k\lesssim \varepsilon^2 n$ as long as $\varepsilon \gg n^{-1/2}$ (see \S 4).

\subsection{Concentration vs Chebyshev's inequality}

The function $c(\varepsilon)\simeq |\log \varepsilon|$ that 
appears in this simple argument of Theorem \ref{thm:w-dvo} is a nice feature. In 
particular, if it is combined with the fact that there exist $n$-dimensional normed spaces $X$ with critical dimension $k(X)\simeq \log n$
whereas $\beta_{\gamma_n}(X) \lesssim n^{-\alpha}$ for some absolute constant $\alpha > 0$ (see Proposition \ref{prop:beta-ell-p}), then 
Theorem \ref{thm:w-dvo} yields the existence of almost Euclidean subspaces of the same dimension 
as Milman's formulation \cite{Mil1} provides, but with better dependence on $\varepsilon$. Moreover, for those 
spaces we conclude that the $k$-dimensional {\it random} subspace is $(1+\varepsilon)$-Euclidean 
with probability $>1-e^{-ck}$ as long as $k\ls c\log n/ \log\frac{1}{\varepsilon}$. 

In the light of the above comments, the random version of Dvoretzky's theorem can be complemented in the following way:

\begin{corollary}
Let $X$ be an $n$-dimensional normed space. Then, for any $\varepsilon \in (0,1)$ there exists $k \gr c \max \left\{ \varepsilon^2 k(X) , (\log \frac{1}{\varepsilon})^{-1} \log \frac{1}{\beta(X)}\right\}$ such that the random $k$-dimensional subspace $F$ satisfies: 
\begin{align*}
\frac{1-\varepsilon}{M} B_F \subseteq B_X\cap F \subseteq \frac{1+\varepsilon}{M} B_F,
\end{align*} with probability greater than $1-e^{-ck}$.
\end{corollary}

\subsection{Two examples on the dependence on $k(X,\varepsilon)$}

The classical results of Milman \cite{Mil1}, Gordon \cite{Go} and Schechtman \cite{Sch1} predict that the dependence on $\varepsilon$ 
should be $\varepsilon^2$. Here we show that this is always the case after a small perturbation. 

Let $f: \mathbb R^{n+1}\to \mathbb R$ be the mapping:
\begin{align}
f(t,x)=  |t|+\|x\|_\infty, \quad (t,x)\in \mathbb R \times \mathbb R^n. 
\end{align} Then, we have the following properties:

\begin{itemize}
\item [i.] $f$ is $2$-equivalent norm to $\|\cdot\|_{\ell_\infty^{n+1} }$.
\item  [ii.] ${\rm Var}[f(Z)]\simeq 1$.
\item [iii.] Moreover, we have: $k(f)\simeq k(\ell_\infty^n)\simeq \log n$.
\item [iv.] $\mathbb E [f(Z)] \simeq \sqrt{\log n}$. 
\end{itemize}

\smallskip

\noindent  We also have the following:

\begin{proposition} \label{prop:ex-1-moms}
For any $r\gr 1$ we have:
\begin{align*}
\left( \mathbb E \left| f(g_1,Z_1)-f(g_2,Z_2) \right|^r \right)^{1/r}\simeq \sqrt{r},
\end{align*} where $Z_1,Z_2$ are independent, standard Gaussian random vectors on $\mathbb R^n$.
\end{proposition}

\noindent {\it Proof.} Using the triangle inequality we may write:
\begin{align*}
\left( \mathbb E \left| f(g_1,Z_1)-f(g_2,Z_2) \right|^r \right)^{1/r} &\gr \left( \mathbb E \left| |g_1|-|g_2| \right|^r\right)^{1/r}- 
\left(\mathbb E\left| \|Z_1\|_\infty-\|Z_2\|_\infty \right|^r \right)^{1/r} \\
&\gr c_2\sqrt{r} - C_2\frac{r}{\sqrt{\log n} } \gr c_2' \sqrt{r}
\end{align*} for all $1 \ls r\ls c_3\log n$. The assertion follows. \prend

\begin{corollary} \label{cor:ex-1-conc}
Let $f$ be as above. Then, we have:
\begin{align*}
\mathbb P \left( \big | f(Z)-\mathbb E f(Z) \big | >\varepsilon \mathbb Ef(Z) \right) \gr c_1e^{-C_1 \varepsilon^2 \log n},
\end{align*} for all $\varepsilon>0$, where $Z\sim N({\bf 0}, I_{n+1})$.
\end{corollary}

\noindent {\it Proof.} It follows from Proposition \ref{prop:ex-1-moms} and the Paley-Zygmund inequality. \prend

\medskip

Arguing as in \cite[Section 5]{PVZ} and using Corollary \ref{cor:ex-1-conc} we may conclude the following:

\begin{theorem} \label{thm:ex1}
For every $n$ there exists an $n$-dimensional, 1-unconditional normed space $X$ which is $2$-isomorphic to $\ell_\infty^n$ 
and has the following property: If $k(X,\varepsilon)$ is the maximal dimension $k$ for which the $k$-dimensional random
 subspace of $X$ is $(1+\varepsilon)$-Euclidean with probability greater than $1-e^{-k}$, then $k(X, \varepsilon)\simeq \varepsilon^2 \log n$.
\end{theorem}

At this point let us note that recently Tikhomirov \cite{Tik-1-unc} showed that there exists 1-unconditional normed space $X$ whose ball
is in John's position and the critical dimension $k(X,\varepsilon)$ in the randomized Dvoretzky is of the order $\varepsilon^2 \log n$.

In the classical paper of Figiel, Lindenstrauss and Milman \cite{FLM}, it is proven that for spaces with cotype $2\ls q<\infty$ (see e.g. \cite{MS}
for the related definition) the critical dimension in the randomized Dvoretzky, in John's position, is at least of the order $\simeq_q \varepsilon^2 n^{2/q}$
(after employing Gordon's result from \cite{Go}). The next example is concerned with the dependence on $\varepsilon$ for 
spaces having the cotype property. It shows that there exists an 1-unconditional $n$-dimensional normed space $X$ with cotype $q,\; 2\ls q <\infty$ which has Dvoretzky number $\simeq_q n^{2/q}$, though the critical dimension $k(X,\varepsilon)$ in the randomized 
Dvoretzky is of the exact order $\simeq_q \varepsilon^2n^{2/q}$.

Let $2<q<\infty$. As before, we consider the map $f: \mathbb R^{n+1} \to \mathbb R$ with 
\begin{align}
f(t,x)= |t| +\|x\|_q, \quad (t,x)\in \mathbb R \times \mathbb R^n.
\end{align} We have the following properties:

\begin{itemize}
\item [i.] $\| (t,z)\|_q \ls f(t,z) \ls 2 \| (t,z)\|_q$ for all $t\in \mathbb R$ and $z\in \mathbb R^n$. Therefore, the $q$-cotype constant $C_q(X)$ 
of $X$ satisfies:
\begin{align*}
C_q(X) \simeq C_q(\ell_q^{n+1}) \simeq \sqrt{q}, 
\end{align*} e.g. see \cite{AK}. 
\item [ii.] If $X =(\mathbb R^{n+1}, f(\cdot) )$, then we have: $k(X) \simeq qn^{2/q}$.
\item [iii.] If $Z$ is a standard Gaussian vector on $\mathbb R^{n+1}$, we have: $\mathbb E [f(Z)] \simeq 
\mathbb E\|Z\|_q \simeq \sqrt{q}n^{1/q}$.
\end{itemize}

Recall the following fact proved in \cite{PVZ}:

\begin{lemma}
Let $2<q<\infty$. Then, for all large enough $n$ we have:
\begin{align*}
\left(\mathbb E \big|  \|Z_1\|_q - \|Z_2\|_q \big |^r \right)^{1/r} \ls c_1(q) \mathbb E \|Z_1\|_q \max \left\{ \sqrt{\frac{r}{n}} , \frac{r^{q/2}}{n} \right\}, 
\end{align*} for all $r\gr 1$, where $c_1(q)>0$ is a constant depending only on $q$ and $Z_1,Z_2$ are independent standard Gaussian random 
vectors on $\mathbb R^n$.
\end{lemma}

Using that we obtain the following:

\begin{proposition}
Let $f$ be as above. Then, we have:
\begin{align*}
\left( \mathbb E \big| f(g_1,Z_1) -f(g_2,Z_2)\big|^r\right)^{1/r} \gr c_2(q) \sqrt{ \frac{r}{k(X)}} \mathbb E [f(g_1,Z_1)],
\end{align*} for all $r\gr 1$, where $g_1,g_2$ are i.i.d. standard normals and $Z_1,Z_2$ are independent 
standard Gaussian vectors on $\mathbb R^n$.
\end{proposition}

The proof is similar to that of Proposition \ref{prop:ex-1-moms}, thus it is omitted. Finally, we get:

\begin{corollary}
Let $f$ as above. Then, we have:
\begin{align*}
\mathbb P \left( |f(Z) -\mathbb E f(Z)| > \varepsilon \mathbb E f(Z) \right) \gr c e^{-c_3(q) \varepsilon^2 n^{2/q} },
\end{align*} for all $\varepsilon>0$, where $Z\sim N({\bf 0},I_{n+1})$.
\end{corollary}

The next Theorem is the analogue of Theorem \ref{thm:ex1} in the case of spaces with cotype:

\begin{theorem}
Let $2<q<\infty$ there exists a constant $C(q)\gg 1$ with the following property: For any $n \gr C(q)$ there exists an $n$-dimensional, 1-unconditional normed space $X$ which satisfies:
\begin{itemize}
\item [i.] $X$ has cotype $q$ with $C_q(X)\simeq \sqrt{q}$ and,
\item [ii.] If for any $\varepsilon\in (0,1)$ we denote by $k(X,\varepsilon)$ the largest dimension $k$ for which the random subspace of $X$ is 
$(1+\varepsilon)$-Euclidean with probabability greater than $1-e^{-ck}$, then $k(X,\varepsilon) \simeq_q \varepsilon^2 n^{2/q}$.
\end{itemize}
\end{theorem}

A more detailed study of the instability of the concentration will appear in \cite{V-tight}.

\medskip
\section{On the parameter $\beta$}

In this Section we study the parameter $\beta$ and we show connections with the one-sided 
inclusion in the randomized Dvoretzky theorem. Let us recall the definition of $\beta$ for a normed 
space $X=(\mathbb R^n,\|\cdot\|)$ equipped with some log-concave probability measure $\mu$:
\begin{align*}
\beta_\mu(X)= \frac{{\rm Var}_\mu \|Z\|}{(\mathbb E_\mu\|Z\|)^2}.
\end{align*} The first result 
shows that the extremal space for the parameter $\beta(X)$ is the Euclidean:

\begin{lemma}
Let $X$ be $n$-dimensional normed space. Then, one has:
\begin{align*}
\beta (X) \gr \beta (\ell_2^n) \simeq 1/n.
\end{align*}
\end{lemma}

\noindent {\it Proof.} Using polar coordinates we see that:
 \begin{align*}
 \mathbb E\|Z\|^p = \mathbb E\|Z\|_2^p \int_{S^{n-1}} \|\theta\|^p \, d\sigma(\theta),
 \end{align*} for all $p>0$. Thus, we may write:
 \begin{align*}
 1+\beta(X) = \frac{\mathbb E \|Z\|^2}{(\mathbb E\|Z\|)^2} &=\frac{\mathbb E\|Z\|_2^2}{( \mathbb E \|Z\|_2)^2} \frac{\int \|\theta\|^2\, d\sigma(\theta)}{(\int \|\theta\| \, d\sigma(\theta))^2} \\
 &=(1+\beta(\ell_2^n)) (1+\beta_\sigma(X)).
 \end{align*} Since $\beta_\sigma \gr 0$, the assertion readily follows. \prend
  
\medskip  
 
The following estimates for the parameter $\beta$ of the classical spaces can be found in \cite{PVZ}:
 
\begin{proposition} \label{prop:beta-ell-p}
There exist absolute constants $0<c<1<C$ with the following property: For all $n\gr 1$ sufficiently large, one has
\begin{align*}
\beta (\ell_p^n) \simeq \left\{ \begin{array}{cc}
\frac{2^p}{p^2 n}, & 1\ls p\ls c\log n \\
\frac{1}{ (\log n)^2}, & C\log n \ls p\ls \infty
\end{array} \right. .
\end{align*}
\end{proposition}

Another class of spaces with this property consists of those which ``sit" inside $L_p$ and $p$ might be moderately growing 
along with the dimension of the underlying space. Namely, if we take into account the estimate for the variance for the subspaces 
of $L_p$, proved in \cite{PV-dvo-L_p}, it follows that for all 
$1\ls p <\infty $ and for all subspaces $X$ of $L_p$ with $\dim X=n$, there exists a linear image $\tilde X$ of $X$ such that: 
\begin{align}
\beta (\tilde X) \ls \frac{C^p}{n} \quad  {\rm or}    \min_{T \in GL(n)} \beta(TX) \ls \frac{C^p}{n} .
\end{align} 

\medskip

General upper and lower bounds for the global parameter $\beta$ are given in the next:

\begin{proposition} \label{prop: bound-k-â}
Let $X$ be finite-dimensional normed space. Then, we have:
\begin{align*}
\log \frac{1}{\beta(X)} \lesssim k(X) \lesssim \frac{1}{\beta(X)}.
\end{align*}
\end{proposition}

Note that the upper bound follows from Claim 2 in \S 3.1. and the definition of $k(X)$. For the lower bound we give
a proof which lies on results of independent interest.

For any norm $\| \cdot \|$ (or any symmetric convex body $K$) on $\mathbb R^n$ we define the following collections of subspaces: 
for $1\ls k \ls n-1$ and $0<\varepsilon, \delta<1$ let
\begin{align}
\mathcal S_{k,\varepsilon} : = \left \{ F\in G_{n,k} \, : \, F \; {\rm is} \; (1+\varepsilon)-{\rm spherical \; subpace \; of} \, X=(\mathbb R^n, \|\cdot\|) \right \}
\end{align} and 
\begin{align}
\mathcal F_{k,\delta} := \left \{ E\in G_{n,k} : (1-\delta) M \ls \|\phi\| \ls (1+\delta) M,  \; \forall \phi\in S_E \right\}.
\end{align}

Next Proposition shows that these two descriptions are essentially equivalent:

\begin{proposition} \label{prop: propo-5-3}
Let $X=(\mathbb R^n, \|\cdot\|)$, let $0<\varepsilon <1$ and let $1\ls m\ls n-1$. Then, we have:
\begin{itemize}
\item [i.] $\mathcal F_{m, \varepsilon/2} \subset \mathcal S_{m, 2\varepsilon}$.
\item [ii.] $\nu_{n,m}(\mathcal S_{m, \varepsilon/3}) \ls C\exp(-c \varepsilon^2 mk(X) )+ \nu_{n,m}( \mathcal F_{m,\varepsilon})$.
\end{itemize}
\end{proposition}

For the proof we shall use a concentration result for the map $F \mapsto M_F$ which is defined as follows: Fix $A$ a symmetric 
convex body on $\mathbb R^n$. Then $F\mapsto M_F(A)=M(A\cap F)$ for any 
$F\in G_{n,k}$, $1\ls k\ls n-1$. Let us mention that a large deviation estimate 
for this mapping was proved by Klartag and Vershynin in \cite[Lemma 3.2]{KV}. Their result reads as follows:

\begin{lemma} \label{lem:KV}
Let $A$ be a symmetric convex body on $\mathbb R^n$ and let $1\ls k\ls n-1$. Then, we have:
\begin{align*}
\nu_{n,k} \left( \left\{ F\in G_{n,k} : M_F(A) > c t M(A) \right\} \right) \ls Ce^{-ct^2k},
\end{align*} for all $t>1$. 
\end{lemma} Note that this estimate, when $k=1$, is reduced to the large deviation estimate of 
the norm $\theta \mapsto \|\theta\|_A$. But in that case the above estimate loses by a term $k(A)$ on the exponent. 
In fact, one can recover this missing factor and moreover, prove a concentration inequality for this mapping:

\begin{theorem} \label{thm: conc-M}
Let $A$ be a symmetric convex body in $\mathbb R^n$ and let $1\ls s\ls n-1$. Then, we have:
\begin{align*}
\nu_{n,s} \left( \{ F \in G_{n,s} : \left| M_F(A)-M(A)\right| >\varepsilon M(A) \right\}) \ls C\exp(-c s \varepsilon^2  k(X)),
\end{align*} for all $\varepsilon>0$ where $C,c>0$ are absolute constants.
\end{theorem}

We believe that distributional inequalities for functionals on the Grassmannian are interesting on their own right, thus we discuss this
topic separately. In fact, we provide two proofs of the latter probabilistic estimate in \S 6. 
Taking this for granted we may proceed with the proof of Proposition \ref{prop: propo-5-3}.

\medskip

\noindent {\it Proof of Proposition \ref{prop: propo-5-3}.} We write $k=k(X)$ for short. The first inclusion is trivial. For the second 
part, note the following inclusion:
\begin{align*}
\mathcal S_{m, \varepsilon} \subset \mathcal A_{m,\varepsilon} :=
\{ E\in G_{n,m} : (1+\varepsilon)^{-1}M_E \ls \|z\| \ls (1+\varepsilon) M_E , \; \forall z\in S_E\}.
\end{align*} If we define
\begin{align*}
\mathcal C_{m,\delta} := \{E\in G_{n,m} : |M_E-M|>\delta M\} , \; 0 < \delta < 1,
\end{align*} then we have:
\begin{align*}
\nu_{n,m}(\mathcal S_{m, \varepsilon/3})  \ls \nu_{n,m} (\mathcal A_{m,\varepsilon/3}) 
&\ls  \nu_{n,m}(\mathcal C_{m, \varepsilon/3}) + \nu_{n,m} (\mathcal A_{m ,\varepsilon/3} \setminus \mathcal C_{m, \varepsilon/3}) \\
&\ls  C e^{-c\varepsilon^2 m k} + \nu_{n,m}( \mathcal F_{m,\varepsilon}),
\end{align*} where we have used Theorem \ref{thm: conc-M} and the fact that 
$\mathcal A_{m, \varepsilon/3} \setminus \mathcal C_{m, \varepsilon/3} \subset \mathcal F_{m, \varepsilon}$. \prend

\bigskip

\noindent {\it Proof of Proposition \ref{prop: bound-k-â}.} (Lower bound). We set $\beta=\beta(X)$. Note that 
Theorem \ref{thm:w-dvo} implies that the set
\begin{align*}
\mathcal S_{s,1/6}= \left \{ F\in G_{n,s} \, : \, F \; {\rm is} \; 7/6-{\rm spherical \; subpace \; of} \, X \right\},
\end{align*} has $\nu_{n,s}( \mathcal S_{s,1/6})\gr 1-e^{-s}$ for $s\simeq |\log\beta|$. 
Furthermore, Proposition \ref{prop: propo-5-3} implies that:
\begin{align*}
1-e^{-c_1s} \ls \nu_{n,s}\left( \left\{ E \in G_{n,s} : \frac{M}{2} \ls \|z\| \ls 2M, \; \forall z\in S_E \right\}\right).
\end{align*} Now we may employ a result of Milman and Schechtman from \cite{MS-global} (see also \cite{HW} for a recent development) to conclude. \prend

\medskip

\noindent {\it Note.} Let us mention that the bounds we derive in Proposition \ref{prop: bound-k-â} are sharp (up to constants): 
For $X = \ell_1^n$ note that we have $ k(\ell_1^n ) \simeq 1/ \beta(\ell_1^n )$. The example of $\ell_p^n$ with $p =c_0 \log n$ for some sufficiently small absolute constant $c_0>0$ (see Proposition \ref{prop:beta-ell-p}) guarantees the existence of
$n$-dimensional normed spaces $X$ with critical dimension $k(X)\simeq \log n$ and $\beta(X)\lesssim n^{-\alpha}$.

\subsection{Dvoretzky's theorem revisited} In \cite{PVsd} we prove a refined random version of the classical dimension reduction lemma
of Johnson and Lindenstrauss \cite{JL} in terms of the parameter $\beta$.
More precisely, given a finite dimensional normed space $X=(\mathbb R^n,\|\cdot\|)$ we may define the 
parameter ${\bf ov}(X)$ associated with
$X$ as follows:
\begin{align*}
{\bf ov}(X)^2:= \beta(X) k(X) = \frac{  {\rm Var}\|Z\|  }{b(X)^2},
\end{align*} where $Z\sim N({\bf 0} ,I_n)$.
Then, the applications in \cite{PVsd} may exhibit improved one-sided behavior which takes into account the parameter ${\bf ov}(X)$. 
It turns out that the parameter ${\bf ov}(X)$ measures the sharpness of the concentration of the norm $\|\cdot\|_X$ in terms of the 
Lipschitz constant; in particular, it measures how ``over-concentrated" is the norm (see \cite{V-tight} for the related subject).

Below we give a more detailed formulation of the random version of Dvoretzky's theorem in terms of this parameter. We also provide
one-sided, almost isometric inclusion in large dimensions proportional to $1/\beta$. The latter can be viewed as the almost isometric version of the result of Klartag and Vershynin in \cite{KV}. For presentation purposes we choose to state and prove the theorem in the dual setting. For this
end we set for any symmetric convex body $A$ in $\mathbb R^n$ the following (dual) parameters $w(A)=M(A^\circ)$, $k_\ast(A)=k(A^\circ)$, 
$\beta_\ast(A)=\beta(A^\circ)$ and ${\bf ov}_\ast(A) = {\bf ov}(A^\circ)$. With this notation we have the following:

\begin{theorem}[Randomized Dvoretzky] \label{thm:new-Dvo}
For any symmetric convex body $A$ we set $k_\ast:=k_\ast(A)$, $\beta_\ast:=\beta_\ast(A)$ and 
$[{\bf ov}_\ast(A)]^2=  \beta_\ast(A)k_\ast(A)$. Then, for every $\varepsilon\in (0,1)$ we have:
\begin{itemize}
\item [\rm i.] There exists a set $\mathcal A$ in $G_{n,k}$ with $k \simeq \varepsilon^2 k_\ast$ and 
$\nu_{n,k}(\mathcal A) \gr 1-e^{-c\varepsilon^2 k_\ast}$
such that
\begin{align*}
P_E(A) \subseteq (1+\varepsilon) w(A) B_E,
\end{align*} for all $E\in \mathcal A$.
\item [\rm ii.] There exists a set $\mathcal B$ in $G_{n,\ell}$ with $\ell \simeq \frac{ \varepsilon^2}{ \beta_\ast \log (1/ \varepsilon)}$ and 
$\nu_{n,\ell}(\mathcal B) \gr 1- e^{-c \varepsilon^2 /\beta_\ast}$ such that
\begin{align*}
P_F(A) \supseteq (1-\varepsilon) w(A) B_F,
\end{align*} for all $F\in \mathcal B$.
\item [\rm iii.] There exists a set $\mathcal C$ in $G_{n,m}$ with $m \simeq \frac{\varepsilon^2 k_\ast}{ \log (1 / \varepsilon) } $ and 
$\nu_{n,m}(\mathcal C)\gr 1-e^{-c\varepsilon^2 k_\ast}$ such that
\begin{align*}
w(A) \left(1 - \frac{ \varepsilon  {\bf ov}_\ast(A) }{\log (\frac{1}{\varepsilon}) } \log \left(\frac{e}{\varepsilon {\bf ov}_\ast(A) } \right) \right)B_G
 \subseteq P_G(A) \subseteq \left(1+\frac{\varepsilon}{\sqrt{\log (1/\varepsilon)}} \right) w(A) B_G,
\end{align*} for all $G\in \mathcal C$,
\end{itemize} where $c>0$ is an absolute constant.
\end{theorem}

We shall need some auxiliary results. The first one follows from Klartag-Vershynin estimate (Lemma \ref{lem:KV}) in the dual setting.

\begin{lemma} \label{lem:stab-w}
Let $A$ be a symmetric convex body on $\mathbb R^n$ and let $1\ls k\ls n-1$. Then, we have:
\begin{align*}
\left( \int_{G_{n,k}} ( w(P_FA))^k \, d\nu_{n,k}(F) \right)^{1/k} \ls C w(A) ,
\end{align*} where $C>0$ is an absolute constant.
\end{lemma}

The next Lemma is essentially from \cite{KV}. However, the next formulation, takes into account the magnitude of the constants involved.

\begin{lemma} [dimension lift] \label{lem:dim-lift}
Let $A$ be a symmetric convex body on $\mathbb R^n$ and let $1\ls s\ls n-1$. Then, for any $p \gr s$ we have:
\begin{align*}
\left( \int_{G_{n,s}} [r(P_EA)]^{-p/2} \, d\nu_{n,s}(E) \right)^{2/p} &\ls 
\frac{ e^{cs/p} a_{s,p} }{w(A)}  \left(\frac{w(A)}{w_{-2p}(A)} \right)^2,
\end{align*} where $c>0$ is an absolute constant and $a_{s,p}$ is defined in Lemma \ref{lem: stability of moms}.
\end{lemma}

\noindent {\it Proof.} Recall that Lemma \ref{lem: stability of moms} yields:
\begin{align*}
J_q(K) =\left(\frac{1}{|K|} \int_K \|x\|_2^q  \, dx \right)^{1/q} \gr a_{n,q}^{-1} R(K), 
\quad a_{n,q}^{-q} = \frac{q}{2}B(q,n+1).
\end{align*} Also Lemma \ref{lem:identity} implies:
\begin{align*}
{\rm vrad}(K)^n J_q^q(K)=M_{-(n+q)}^{-(n+q)}(K).
\end{align*} Recall that $1/r(K)=R(K^\circ)$, hence if we apply the above for $K=(P_EA)^\circ$ we get:
\begin{align*}
\int_{G_{n,s}} \frac{1}{r(P_EA)^{p/2} } \, d\nu_{n,s}(E) &\ls a_{s,p}^{p/2} \int_{G_{n,s}} J_p^{p/2} ((P_EA)^\circ)\, d\nu_{n,s}(E) \\
&= \frac{p a_{s,p}^{p/2} }{s+p} \int_{G_{n,s}} w_{-(s+p)}^{-(s+p)/2}(P_EA) \cdot \left[{\rm vrad}((P_EA)^\circ)\right]^{-s/2} \, d\nu_{n,s}(E) \\
&\ls a_{s,p}^{p/2} \left( \int_{G_{n,s}} w_{-(s+p)}^{-(s+p)}(P_EA)\, d\nu_{n,s}(E) \right)^{1/2} \cdot 
\left( \int_{G_{n,s}} [{\rm vrad}( (P_EA)^\circ )]^{-s} \, d\nu_{n,s}(E) \right)^{1/2} \\
&\ls  a_{s,p}^{p/2}  w_{-(s+p)}^{-(s+p)/2}(A) \left( \int_{G_{n,s}} [w(P_EA)]^s \, d\nu_{n,s}(E) \right)^{1/2} \\
&\ls a_{s,p}^{p/2} w_{-(s+p)}^{-(s+p)/2}(A) w(A)^{s/2} e^{c s},
\end{align*} where we have also used inequality $w(K){\rm vrad}(K^\circ) \gr 1$ (following by H\"older's inequality), 
and Lemma \ref{lem:stab-w}. The result follows. \prend

\begin{proposition} \label{prop:stab-neg-moms}
Let $A$ be a symmetric convex body on $\mathbb R^n$. Then, we have:
\begin{align*}
\frac{w(A)}{w_{-q}(A)} \ls \frac{ \mathbb E [h_A(Z)]}{( \mathbb E [h_A^{-q}(Z)] )^{-1/q}} \ls 
\exp \left( C \max \{  \sqrt{\beta_\ast} ,q \beta_\ast \} \right),
\end{align*} for all $0<q < c/\beta_\ast$, where $\beta_\ast= \beta_\ast(A)$ and $Z\sim N({\bf 0},I_n)$.
\end{proposition}

\medskip

For a proof of the estimate on the Gaussian averages the reader may consult \cite[Corollary 3.5]{PVsd}. In order to pass to the sphere, recall the well known formula (which follows by integration in polar coordinates):
\begin{align*}
\left( \mathbb E [h_A^r(Z)] \right)^{1/r}= w_r(A) \left( \mathbb E \|Z\|_2^r\right)^{1/r}, \quad Z\sim N({\bf 0},I_n),
\end{align*} for all $-n<r\neq 0$. Thus, H\"older's inequality yields the asserted bound. 

\medskip

Now we turn to proving Theorem \ref{thm:new-Dvo}.

\smallskip

\noindent {\it Proof of Theorem \ref{thm:new-Dvo}.} Fix $1\ls s \ls n-1$. Using Lemma \ref{lem:dim-lift} and 
restricting ourselves to the range $1\ls s\ls p\ls c_1/\beta_\ast$ we obtain:
\begin{align*}
\left( \int_{G_{n,s}} [r(P_EA)]^{-p} \, d\nu_{n,s}(E) \right)^{1/p} &
\ls \frac{1}{w(A)} \exp\left( \frac{cs}{p} \log \left(\frac{ep}{s} \right) + c\sqrt{\beta_\ast} + cp\beta_\ast \right),
\end{align*} where we have also used Proposition \ref{prop:stab-neg-moms}. 

\smallskip

\noindent (ii) Given $\varepsilon\in (0,1)$ with $\varepsilon \gtrsim \sqrt{\beta_\ast}$ (otherwise there is nothing to prove) we 
choose $p\simeq \frac{\varepsilon}{\beta_\ast} $ and $s\simeq \varepsilon p / \log(1/\varepsilon)$. Then, we obtain:
 \begin{align*}
 \left( \int_{G_{n,\ell}} [r(P_EA)]^{-p} \, d\nu_{n,\ell}(E) \right)^{1/p} &\ls 
\frac{1}{w(A)} \exp \left( C\varepsilon \right),
\end{align*} with $\ell =s \simeq \varepsilon^2/ ( \beta_\ast\log \frac{1}{\varepsilon})$. The asserted set $\mathcal B$ is defined as:
\begin{align*}
\mathcal B:= \left\{ F \in G_{n,\ell} : r(P_FA) \gr w(A) e^{-2C\varepsilon} \right\}.
\end{align*} Note that by Markov's inequality we have: $\nu_{n,\ell}(\mathcal B^c)\ls e^{-c_1\varepsilon p}\ls e^{-c_2\varepsilon^2 /\beta_\ast}$
 
\smallskip
 
\noindent  (iii) Set $\tau_\ast=[{\bf ov}_\ast(A)]^2$. Given $\varepsilon\in (0,1)$ (with $\varepsilon \gtrsim 1/\sqrt{k_\ast}$) we choose 
$p \simeq \varepsilon \sqrt{\tau_\ast} /\beta_\ast $ and 
$s\simeq p\varepsilon \sqrt{\tau_\ast} / \log \frac{1}{\varepsilon} $ and we argue as before to get:
\begin{align*}
\mathcal C:= \left\{ G\in G_{n,m} : r(P_G(A)) \gr w(A) e^{-2C\varepsilon \sqrt{\tau_\ast} \log(\frac{e}{\varepsilon \tau_\ast}) /\log \frac{1}{\varepsilon} } \right\},
\end{align*} where $m =s \simeq \varepsilon^2k_\ast / \log\frac{1}{\varepsilon}$ and 
$\nu_{n,m}(\mathcal C^c) \ls e^{-c_1 p\frac{\varepsilon\sqrt{\tau_\ast}}{\log(1/\varepsilon)} \log (\frac{e}{\varepsilon \sqrt{\tau_\ast}}) } \ls 
e^{-c_2 \varepsilon^2 k_\ast}$. The proof is complete. \prend

\medskip
\section{Small deviation estimates and the parameter $\vartheta$}

The aim of this paragraph is to prove tight small deviation and small ball probability estimates 
and illustrate the critical parameters that govern such estimates. The technique can be applied to the 
more general context of the log-concave probability measures. Our starting point is 
the following one-sided small deviation inequality proved in \cite{PVsd}:

\begin{theorem} \label{thm:sd-gauss}
Let $f:\mathbb R^n\to \mathbb R$ be a convex map with $f\in L_1(\gamma_n)$. Then, 
\begin{align} \label{eq:gsd}
\mathbb P \left( f(Z) < m-t \right)\ls \frac{1}{2}\exp(-ct^2/ \|f-m\|_{L_1}^2),
\end{align} for all $t>0$, where $m$ is a median of $f(Z)$ and $Z\sim N({\bf 0}, I_n)$.
\end{theorem} 

The latter can be viewed as a strengthening of the one-sided classical Gaussian concentration, inside the class of convex functions,
since one can replace the Lipschitz constant by the a priori smaller quantity of the standard deviation. 
The main ingredients of the proof can be summarized in the following two observations:

\begin{itemize}
\item When $f$ is a convex function, the map $t\mapsto \Phi^{-1} \circ \gamma_n (f \ls t)$ is concave, 
(this follows by Ehrhard's inequality \cite{Ehr}) and,
\item the derivative of the above function at the median $m$ of $f$ with respect to $\gamma_n$ is relatively large in terms of the
standard deviation of $f$.
\end{itemize}
Here we extend the above approach in the context of log-concave probability measures. 
Let $\mu$ be a log-concave, Borel probability measure on $\mathbb R^n$. Let $Z$ be a random vector 
on $\mathbb R^n$ distributed according to $\mu$, i.e. $\mathbb P(Z \in B)=\mu(B)$ for any Borel set $B\subseteq \mathbb R^n$.
 
For any convex function $\psi: \mathbb R^n\to \mathbb R$ we define the random variable $\xi:=\psi(Z)$. We write
$F_\xi(t): =\mu( z\in \mathbb R^n : \psi(z) \ls t) = \mathbb P( \psi(Z) \ls t)$ for the cumulative distribution function of $\xi$. The first main step
is to replace the function $\Phi$ by a suitably chosen function that yields the concavity. The next
standard Lemma (for a proof e.g. see \cite{Bob}) shows that in the context of log-concave measures this function can be 
at least the exponential.

\begin{lemma} \label{lem:cdf}
Let $\mu$ be a log-concave probability measure on $\mathbb R^n$ and let $\psi$ be a convex map on $\mathbb R^n$. 
If $F(t)= \mathbb P( \psi(Z) \ls t)$, then we have the following:
\begin{itemize}
\item [a.] For $t,s\in \mathbb R$ and $0<\lambda<1$ we have: 
\begin{align*}
F \left( (1-\lambda)t +\lambda s \right) \gr \left[ F(t)\right]^{1-\lambda} \left[F(s) \right]^\lambda,
\end{align*} that is $F$ is log-concave.

\item [b.] If $\psi$ is a semi-norm\footnote{ A seminorm $\psi :V\to [0,\infty)$ on a vector space $V$ is a function which is positively homogeneous, i.e. $\psi(\lambda v)=|\lambda| \psi(v)$ for all scalars $\lambda$ and $v\in V$ and sub-additive, that is $\psi(u+v) \ls \psi(u) +\psi(v)$
for all $u,v\in V$.}, then we also have:
\begin{align*}
1-F\left( (1-\lambda)t -\lambda s\right) \gr (1-F(t))^{1-\lambda} F(s)^\lambda
\end{align*} for all $t,s>0$ and $0<\lambda < \frac{t}{t+s}$.
\end{itemize}
\end{lemma}

\noindent {\it Sketch of proof.} For b. we may check (using the subadditivity of the seminorm) that 
\begin{align*}
(1-\lambda) \{\psi > t\} + \lambda \{ \psi \ls s\} \subset  \{ \psi > (1-\lambda) t -\lambda s\},
\end{align*} and we use the log-concavity of $\mu$. \prend

\smallskip

\begin{remark} Note that (b) easily implies Borell's lemma from \cite{Bor}: Fix $s>0$ (say $s=1$). For any $t>s$ we choose $\lambda\in (0,\frac{t}{t+s})$ such that $(1-\lambda)t-\lambda s=s$, i.e. $\lambda= \frac{t-s}{t+s}$ and hence,
\begin{align*}
1-F(t) \ls F(s) \left( \frac{1-F(s)}{F(s)} \right)^{\frac{t+s}{2s}}, \quad 0<s<t.
\end{align*}
\end{remark}

\medskip

The next lemma shows that $f(m)$ can be estimated in terms of the standard deviation of the function, thus fulfills the second 
main observation in the general context of log-concave measures. 

\begin{lemma} \label{lem:der-var} 
Let $\mu$ be a log-concave probability measure on $\mathbb R^n$ and let $\psi : \mathbb R^n \to \mathbb R$ be a convex map 
with $\psi \in L_1(\mu)$. 
If $F(t)= \mathbb P( \psi(Z) \ls t), \, t\in \mathbb R$ and $f=F'$, then
\begin{align*}
f(m) \gr \frac{1}{32 \| (\psi-m)_+ \|_{L_1(\mu)} },
\end{align*} where $m$ is a median of $\psi$ with respect to $\mu$.
\end{lemma}

\noindent {\it Proof.} Since $F$ is log-concave (Lemma \ref{lem:cdf}) we have:
\begin{align*}
2f(m)=(\log F)'(m) &\gr \frac{\log F(m+\delta)-\log F(m)}{\delta} \\
&= \frac{1}{\delta} \log \left[ 1 + 2 \mathbb P\left( m< \psi(Z) \ls m+\delta \right) \right] \\
&\gr \frac{1}{\delta} \mathbb P\left( m< \psi(Z) \ls m+\delta \right)\\
&= \frac{1}{\delta} \left[\frac{1}{2} - \mathbb P(\psi(Z) >m+\delta) \right].
\end{align*} The choice $\delta= 4 \| (\psi -m)_+ \|_{L_1}$ and Markov's inequality yield the result.\prend

\begin{note}
In the above argument we may replace the $L_1$ norm by the difference $m_q-m$ for any $q$-quantile $m_q$ of 
$\psi(Z)$ with $q>1/2$. That is, if $\mathbb P(\psi(Z) \ls m_q)=q>1/2$ then $f(m) \gr \frac{q-1/2}{m_q-m}$. 
Note that if $\psi \in L_{p,\infty}(\mu)$ then $m_q-m \ls (1-q)^{-1/p} \|(\psi-m)_+\|_{p,\infty}$.
\end{note}

\bigskip

\begin{theorem} \label{thm:main-2}
Let $\mu$ be a log-concave probability measure on $\mathbb R^n$. Let $\psi$ be a convex function on $\mathbb R^n$ 
and let $m$ be a median for $\psi$ with respect to $\mu$. Then, 
\begin{align*}
\mu \left( \left \{ x: \psi(x) < m -t \int |\psi-m| \, d\mu \right \} \right) \ls \frac{1}{2} \exp(-t/16),
\end{align*} for all $t>0$.
\end{theorem}

\noindent {\it Proof.} Let $Z$ be a random vector distributed according to $\mu$ and let 
$F(t) = \mathbb P ( \psi(Z) \ls t), \; t\in \mathbb R$ and $f=F'$. Then, by Lemma \ref{lem:cdf} we have:
\begin{align*}
\log F(m-t) - \log F(m) \ls -t\frac{f(m)}{F(m)}  = -2t f(m).
\end{align*} It follows that:
\begin{align} \label{eq:sd-theta}
F(m-t)\ls F(m) \exp(- 2t f(m)) =\frac{1}{2}\exp(- 2t f(m)) ,
\end{align} for all $t>0$. Using Lemma \ref{lem:der-var} we get:
\begin{align*}
F(t-m) \ls \frac{1}{2} \exp\left( -t/ (16 \mathbb E_\mu |\psi(Z)-m| ) \right),
\end{align*} as required. \prend

\bigskip

For any symmetric convex body $K$ in $\mathbb R^n$ and for any log-concave probability measure $\mu$ on $\mathbb R^n$ we
write $F_K$ for the cumulative distribution function of $\|Z\|_K$, where $Z$ is distributed according to $\mu$, i.e. 
$F_K(t)=\mu(x: \|x\|_K\ls t), \; t>0$. If $f_K$ is the density of this random variable, we define:
\begin{align}
\vartheta(\mu, K) \equiv \vartheta_\mu(K) := mf_K(m),
\end{align} where $m$ is a median for $\|Z\|_K$. With this notation Lemma \ref{lem:der-var} shows that 
\begin{align} \label{eq:theta-beta}
\vartheta_\mu(K) \gtrsim \frac{1}{\sqrt{\beta_\mu(K)}}.
\end{align} Below, we also show that $\vartheta$ raises naturally in small deviation
and small ball estimates for log-concave probability measures on $\mathbb R^n$. For this end, let us recall the {\it $B$-property}.
A pair $(\mu,K)$ of a log-concave probability measure on $\mathbb R^n$ and a symmetric convex body $K$ in $\mathbb R^n$ it is said 
to have the $B$-property if the map $t\mapsto \mu(e^tK)$ is log-concave. 
In \cite{CFM}, Cordero-Erausquin, Fradelizi and Maurey proved that $(\gamma_n,K)$ has the $B$-property for any symmetric convex body $K$ 
in $\mathbb R^n$. They also proved that any pair $(\mu,K)$ of an 1-unconditional log-concave 
measure and 1-unconditional convex body on $\mathbb R^n$ also has the $B$-property. In view of Theorem \ref{thm:main-2},
we have the following analogue of \cite[Proposition 3.2]{PVsd} for log-concave measures:

\begin{theorem} \label{thm:sbp-lc}
Let $K$ be a symmetric convex body in $\mathbb R^n$ and let $\mu$ be a log-concave probability measure. If $m$ is the median
of $x\mapsto \|x\|_K$ with respect to $\mu$, then we have:
\begin{itemize}
\item [a.] For every $\varepsilon\in (0,1)$, 
\begin{align*}
\mu(\{x : \|x\|_K \ls (1-\varepsilon)m \} )\ls \frac{1}{2}\exp \left ( -2\varepsilon \vartheta_\mu (K) \right) 
\ls \frac{1}{2} \exp \left(- c\varepsilon / \sqrt{\beta_\mu(K)} \right).
\end{align*}

\item [b.] Furthermore, if the pair $(\mu,K)$ possesses the B-property and $\mu(K)\ls 1/2$, 
\begin{align*}
\mu( \varepsilon K)\ls \varepsilon^{ 2f_K(m) } \mu(K) \ls  \varepsilon^{ \frac{1}{16 \mathbb E | \|Z\|_K-m | }} \mu(K),
\end{align*} for all $\varepsilon\in (0,1)$. In particular,
\begin{align*}
\mu( \{x : \|x\|_K \ls \varepsilon m\})\ls \frac{1}{2} \varepsilon^{2\vartheta_\mu(K)} \ls \frac{1}{2} \varepsilon^{c/\sqrt{\beta_\mu(K)}},
\end{align*} for every $\varepsilon\in (0,1)$.

\end{itemize}
\end{theorem}

\noindent {\it Proof.} The first assertion follows by \eqref{eq:sd-theta} and \eqref{eq:theta-beta}. For the second estimate recall 
the fact that the $B$-property implies that for any symmetric convex body $A$ the map $t\mapsto [\mu(tA)/ \mu (A)]^{1/\log (1/t) }, \; t\in (0,1)$ 
is non-decreasing (e.g. see \cite{LO}). Therefore, if we define
\begin{align*}
d:= \sup_{0<\varepsilon <1} \frac{\log [\mu( \varepsilon K) / \mu(K) ]}{\log (1/\varepsilon) },
\end{align*} the monotonicity shows that:
\begin{align*}
d= \lim_{\varepsilon \to 1^-} \frac{\log [\mu( \varepsilon K) / \mu(K) ]}{\log (1/ \varepsilon) } = -\frac{d}{d\varepsilon} \Big|_{\varepsilon =1} \log F_K(\varepsilon ) = - (\log F_K)'(1)
\end{align*} On the other hand, the map $t\mapsto (\log F_K)'(t)$ is non-increasing (Lemma \ref{lem:cdf}) and since $\mu(K)\ls 1/2$ we get that $m\gr 1$,
which implies
\begin{align*}
(\log F_K)'(1) \gr (\log F_K)'(m) =2f_K(m)
\end{align*}
Combining the above we infer that:
\begin{align*}
\frac{\log [\mu( \varepsilon K) / \mu(K) ]}{\log (1/ \varepsilon) } \ls -2 f_K(m)  \Longrightarrow \mu(\varepsilon K) 
\ls \varepsilon^{2 f_K(m)} \mu(K) 
\end{align*} for every $\varepsilon\in (0,1)$ as required. \prend

\medskip

\begin{remark} \label{rem:optimality-beta}
Let us mention that the small deviation inequality we obtain in terms of the parameter $\beta$ is sharp. For this end consider any symmetric
convex body $K$ in $\mathbb R^n$ and let $\mu_K$ the uniform probability measure concentrated on $K$,
i.e. 
\begin{align*}
\mu_K(A)= |A\cap K|/ |K|, \quad A \subseteq \mathbb R^n, \quad \textrm{Borel}.
\end{align*} If $m$ is the median of $\|\cdot \|_K$ with respect to $\mu_K$, then $|mK|= |K|/2$. Thus, 
\begin{align*}
|\{x\in K : \|x\|_K \ls (1-\varepsilon)m \}|= (1-\varepsilon)^n |mK| = (1-\varepsilon)^n \frac{|K|}{2}.
\end{align*} It follows that,
\begin{align*}
\mu_K(\{x: \|x\|_K \ls (1-\varepsilon)m \}) =\frac{1}{2} (1-\varepsilon)^n \gr \frac{1}{2} e^{-2\varepsilon n}, \quad 0<\varepsilon <1/2.
\end{align*} Finally, note that 
\begin{align*}
\beta_{\mu_K} (\| \cdot \|_K) =\frac{ \int \|x\|_K^2  \, d\mu_K(x) }{ \left(\int \|x\|_K \, d\mu_K(x) \right)^2} - 1 =  
\frac{\frac{n}{n+2}}{\left( \frac{n}{n+1}\right)^2}-1 = \frac{(n+1)^2}{n(n+2)} -1 = \frac{1}{n(n+2)}.
\end{align*} For the latter identities the reader is referred to \cite[\S 2.1]{MiPa}.
\end{remark}

\medskip

\begin{remark}  Recently, in \cite{ENT}, Eskenazis, Nayar and Tkocz proved that the pair $(\nu_1^n,K)$ has the $B$-property for all 
symmetric convex bodies $K$ on $\mathbb R^n$; recall that the measure $\nu_1^n$ has density $d\nu_1^n(x) = 2^{-n} e^{-\|x\|_1}\, dx$.
In view of Theorem \ref{thm:sbp-lc} and their result one has:
\begin{align} \label{eq:exp-sb} 
\nu_1^n(\varepsilon K)\ls  \varepsilon^{c/\mathbb E | \|Z\|_K-m | } \nu_1^n(K), \quad 0<\varepsilon<1,
\end{align} for all symmetric convex bodies $K$ on $\mathbb R^n$ with $\nu_1^n(K)\ls 1/2$. A different estimate was proved in
\cite[Corollary 13]{ENT} in terms of the inradius $r(K)$ of the body $K$, which is in the spirit of \cite{LO}. 
Let us mention that for the class of 1-unconditional convex bodies a better estimate than \eqref{eq:exp-sb} is 
known for the exponential distribution; see \cite[Proposition 3.4]{PVsd}. Furthermore, for the Gaussian distribution $\gamma_n$ and all symmetric convex bodies $K$ in $\mathbb R^n$ with $\gamma_n(K) \ls 1/2$, one has:
\begin{align*}
\gamma_n(tK) \ls (2t)^{cm/(\mathbb E | \|Z\|_K-m|)^2} \gamma_n(K),
\end{align*} for every $t\in [0,1]$, where $m$ is the median of $\|Z\|_K$ with $Z\sim N({\bf 0}, I_n)$ (see \cite{PVsd}). One should compare the latter
with \cite[Theorem 2]{LO}.

In order to illustrate the difference, between estimate \eqref{eq:exp-sb} and the one proved in \cite[Corollary 13]{ENT}, let us consider the following example. Fix $1\ls p \ls \infty$ and let $W$ be an $n$-dimensional exponential random vector, i.e. $W\sim \nu_1^n$. Then, we have
\begin{align*}
{\rm Var}\|W\|_p \simeq_p \left\{ \begin{array}{ll} n^{\frac{2}{p}-1}, & 1\ls p<\infty \\
1 , & p=\infty 
\end{array}\right. ,
\end{align*} for all sufficiently large $n$ (for a proof of this fact the reader may consult \cite[Theorem 3.2 \& \S 3.2]{PVZ}). Furthermore,
\begin{align*}
\mathbb E\|W\|_p \simeq \left\{ \begin{array}{ll} 
pn^{1/p}, & 1\ls p\ls \log n\\
\log n, & p\gr \log n
\end{array} \right. .
\end{align*} For a proof of the latter fact the reader is referred to \cite{SchZ}. Thus, one has:
\begin{align*}
\beta_{\nu_1^n}(B_p^n) \simeq_p  \left\{ \begin{array}{ll} n^{-1}, & 1\ls p<\infty \\
(\log n)^{-2}, & p=\infty
\end{array} \right. ,
\end{align*} for fixed $p$ and for all sufficiently large $n$. It follows that if $m=m_{p,n}$ is the median of $\|W\|_p$ with $1\ls p\ls \infty$, then 
the small ball estimate \eqref{eq:exp-sb} (applied for $K=mB_p^n$) yields 
\begin{align} \label{eq:sb-exp-p}
\mathbb P( \|W\|_p \ls \varepsilon m) \ls \frac{1}{2} \varepsilon^{c/ \sqrt{ \beta_{\nu_1^n} (B_p^n) } },
\end{align} while the corresponding small ball in terms of the inradius $r(B_p^n)$ yields
\begin{align*}
\mathbb P \left( \|W\|_p \ls \varepsilon m\right)\ls \frac{1}{2} \varepsilon^{c \sqrt{k_{\nu_1^n}(B_p^n)}},
\end{align*} where 
\begin{align*}
k_{\nu_1^n}(B_p^n) := \left( r(B_p^n) \cdot \mathbb E\|W\|_p \right)^2\simeq \left\{ \begin{array}{lll} 
n, &  1\ls p\ls 2 \\
p^2n^{2/p}, & 2 \ls p \ls \log n \\
(\log n)^2, & p\gr \log n 
\end{array}\right. .
\end{align*} Actually in view of \cite[Proposition 3.4]{PVsd} one can have an even better estimate than \eqref{eq:sb-exp-p} 
since the $\ell_p$ norms are 1-unconditional.

For further comparison purposes,
recall that if $\mu$ satisfies a Poincar\'e (or spectral gap) inequality with constant $\lambda_1=\lambda_1(\mu)>0$, i.e.
\begin{align*}
\lambda_1 {\rm Var_\mu}(f) \ls \int \|\nabla f\|_2^2 \, d\mu,
\end{align*} for all smooth $f$, then 
\begin{align*}
\sqrt{\lambda_1} \|f-m\|_{L_1(\mu)} \ls \left( \mathbb E_\mu \|\nabla f\|_2^2 \right)^{1/2} \ls \|f\|_{\rm Lip},
\end{align*} for every Lipschitz map $f$, where $m$ is a median of $f$. In the case that $f$ is the norm $\|\cdot\|_K$ induced 
by the symmetric convex body $K$, note that $\|f\|_{\rm Lip}=1/r(K)$, which indicates that the $L_1$ deviation is in general smaller than the inradius of the body. 

\end{remark}

\medskip

\begin{examples} 1. (The parameter $\vartheta_{\gamma_n}(B_\infty^n)$). Let $m$ be the median for $x\mapsto \|x\|_\infty$ 
with respect to $\gamma_n$. Then, we have:
\begin{align*}
\vartheta(\gamma_n, B_\infty^n) \simeq \log n.
\end{align*}
Indeed; note that for any $s>0$ we have:
\begin{align*}
F_{B_\infty^n}(s) = \gamma_n(sB_\infty^n) = (\gamma_1([-s,s]))^n =( 2\Phi(s)-1)^n.
\end{align*} Hence, 
\begin{align*}
f_{B_\infty^n} (s) =2n \phi(s) (2\Phi(s)-1)^{n-1}.
\end{align*} In particular, if $m$ is the median of $x\mapsto  \|x\|_\infty$ with respect to $\gamma_n$, then 
\begin{align*}
f_{B_\infty^n} (m) = n \frac{ \phi(m)}{ 2\Phi (m)-1} = n2^{1/n} \phi( m).
\end{align*}
On the other hand we have:
\begin{align*}
2^{-1/n} =\int_{-m}^m \phi(t)\, dt =1-2 \int_m^{\infty} \phi(t)\, dt \Longrightarrow  2\int_m^\infty \phi =1-2^{-1/n} \simeq 1/n
\end{align*} and from standard estimates on the error function we have:
\begin{align*}
\int_m^\infty \phi(t)\, dt \simeq \frac{1}{m} \phi(m), \quad m\to \infty.
\end{align*} Hence, we obtain $f_{B_\infty^n}(m) \simeq m$ or $\vartheta(\gamma_n, B_\infty^n) \simeq m^2 \simeq  \log n.$

\medskip

\noindent 2. (The parameter $\vartheta_{\gamma_n} (B_2^n)$). Let $m$ be a median for $x \mapsto \|x\|_2$ with respect to Gaussian. 
Then, we have:
\begin{align*}
\vartheta=\vartheta(\gamma_n, B_2^n)\simeq \sqrt{n}.
\end{align*} Indeed; it is known that $\beta(B_2^n)\simeq 1/n$ which implies 
that $\vartheta \gtrsim \sqrt{n}$. For the upper estimate we argue as follows:
\begin{align*}
e^{-c_1\vartheta^2} \gr \gamma_n \left( \frac{m}{2}B_2^n \right) \gr (2\pi)^{-n/2} \left| \frac{m}{2}B_2^n \right| e^{-m^2/8}\gr e^{-c_2n}, 
\end{align*} where we have used the fact that $m\simeq \sqrt{n}$.
\end{examples}

\bigskip

In fact the upper estimate holds for any centrally symmetric convex body in $\mathbb R^n$:

\begin{corollary}
Let $K$ be a symmetric convex body on $\mathbb R^n$. Then,
one has:
\begin{align*} \vartheta(\gamma_n, K) \lesssim \sqrt{n}.
\end{align*}
\end{corollary}

\noindent {\it Proof.} We shall need \cite[Lemma 2.1]{KV} which states:
\begin{align*}
\frac{1}{2}\sigma(S^{n-1} \cap \frac{1}{2}L) \ls \gamma_n(\sqrt{n}L),
\end{align*} for any centrally symmetric convex body $L$. Using the dual Sudakov inequality (see e.g. \cite{LT}) we get:
\begin{align*}
\frac{|\sqrt{n}B_2^n|}{|\sqrt{n}B_2^n \cap \frac{m}{4}K|} \ls N(\sqrt{n}B_2^n, \frac{m}{4}K) \ls \exp\left( c_1(\mathbb E\|Z\|_K)^2 (4\sqrt{n}/m)^2 \right)
\ls e^{c_2n},
\end{align*} where $m$ is a median for $\|Z\|_K$ and $Z\sim N({\bf 0},I_n)$. On the other hand we have:
\begin{align*}
\frac{|\sqrt{n}B_2^n \cap \frac{m}{4}K|}{|\sqrt{n}B_2^n|} =\sigma(S^{n-1}\cap \frac{m}{4\sqrt{n}}K) \ls 2\gamma_n(\frac{m}{2}K) \ls \exp(-c_3\vartheta^2).
\end{align*} Combining all the above we get the result. \prend

\subsection{A remark on isoperimetry}

Let $\nu_1$ be the 1-dimensional exponential measure with density $d\nu_1(x) =\frac{1}{2}e^{-|x|}\, dx$. Let $F_{\nu_1}$ be its cumulative distribution function, i.e.
\begin{align*}
F_{\nu_1}(x) = \nu_1 ( (-\infty, x ]) = \frac{1}{2}\int_{-\infty}^x e^{-|t|}\, dt, \quad x\in \mathbb R.
\end{align*} Let $X=(\mathbb R^n, \|\cdot\|)$ be a normed space and let $\mu$ be a log-concave measure on $\mathbb R^n$. 
Following \cite{Bob} we shall denote the induced measure by $\mu_X$. That is the push forward of the measure $\mu$ under the 
map $x\mapsto \|x\|$, i.e. 
$$\mu_X(A)=\mu( x\in \mathbb R^n : \|x\| \in A), \quad A\subseteq [0,\infty), \; {\rm Borel}.$$

Let $T:(0,\infty)\to \mathbb R$ with $T= F_{\nu_1}^{-1}\circ F_K$ be the map which transports $\mu_X$ to $\nu_1$, where $K=\{x: \|x\| \ls 1\}$. 
Bobkov in \cite{Bob} proved the following:

\begin{proposition} \label{prop:Bob}
Let $\mu$ be a log-concave probability measure on $X=(\mathbb R^n, \|\cdot\|)$. Then, one has:
\begin{align*}
{\rm Is}(\mu_X)=\inf_{s>0} (F_{\nu_1}^{-1}\circ F_K)'(s) \gtrsim \frac{1}{ \mathbb E_\mu \|Z\|}.
\end{align*}
\end{proposition}

It is known that ${\rm Is}(\mu_X) \lesssim 1/\sqrt{ {\rm Var}_\mu\|Z\| }$ (see e.g. \cite{Bob} and the references therein). In \cite[\S 4]{Bob}, Bobkov shows that for $1$-dimensional log-concave measures the reverse estimate also holds true and asks if a reverse estimate 
should be expected for the measure $\mu_X$. Note that even though the measure $\mu_X$ is not necessarily log-concave, 
it enjoys many properties, see e.g. Proposition \ref{prop:Bob}. Here we observe that Lemma \ref{lem:der-var} partially answers the aforementioned question:

\begin{corollary}
Let $\mu$ be a log-concave probability measure on $\mathbb R^n$ and let $K$ be a symmetric convex body on $\mathbb R^n$. 
Then, we have:
\begin{align*}
\inf_{0<s \ls m} (F_{\nu_1}^{-1}\circ F_K)'(s) =(F_{\nu_1}^{-1}\circ F_K)'(m) \gtrsim \frac{1}{\mathbb E_\mu \big| \|Z\|_K - m \big|},
\end{align*} where $m$ is the median of the function $x\mapsto \|x\|_K$ with respect to $\mu$.
\end{corollary}

\bigskip

Moreover, a lower estimate for the full range of $s$ shouldn't be expected as the next argument shows: 
If $\inf_{s>0} T'(s) := 1/L$, then $T:(0,\infty) \to \mathbb R$ is homeomorphism and 
$1/L$-expansive, i.e. $|Tu-Tv| \gr |u-v|/L$ for all $u,v>0$. Thus, the map
$T^{-1}:\mathbb R \to (0,\infty)$ is $L$-Lipschitz and transports the measure $\nu_1$ to $\mu_X$, i.e. for any Borel set 
$A \subseteq \mathbb R$ we have:
\begin{align} \label{eq:transport}
\nu_1(A)= \mu_X(T^{-1}(A)) = \mu \left( x: \|x\|_K \in T^{-1}(A) \right).
\end{align} 

\noindent {\it Claim.} For any $t>0$ one has\footnote{Indeed; if $|z|\ls t$ then, since $T(m)=0$, we have: 
$|T^{-1}(z)-m| = |T^{-1}(z)-T^{-1}(0)|\ls L |z| \ls t L$.}:
\begin{align*}
T^{-1} ( [-t,t]) \subseteq [m-tL ,m+tL].
\end{align*}

In particular, for $A= (-t, t)$, $t>0$ in \eqref{eq:transport} we have: 
\begin{align*}
1- \mathbb P(  | \|Z\|_K -m | > tL ) &= \mu (x: \|x\|_K \in [m-tL, m+tL] ) \\
& \gr \mu \left( x : \|x\|_K \in T^{-1}( [-t, t] ) \right) = \nu_1([-t,t]).
\end{align*} It follows that $\mathbb P( |  \|Z\|_K -m | >tL ) \ls e^{-t}$ for all $t>0$, or equivalently
\begin{align*}
\mathbb P \left( \big| \|Z\|_K-m \big| > \varepsilon m \right) \ls \exp\left( -\varepsilon m/L \right),
\end{align*} for all $\varepsilon >0$. Hence, if $L\simeq \sqrt{ {\rm Var}_\mu \|Z\| }$, we would obtain:
\begin{align*}
\mathbb P \left( \big| \|Z\|_K-m \big| > \varepsilon m \right) \ls \exp\left( - c \varepsilon/ \sqrt{\beta_\mu(K)} \right),
\end{align*} for all $\varepsilon>0$.
Note that for any fixed $2<p<\infty$ one has:
\begin{align*}
c\exp(-C_p \min\{ \varepsilon^2 n, (\varepsilon n)^{2/p} \}) \ls \mathbb P\left( \big| \|Z\|_p-M_{p,n} \big|> \varepsilon M_{p,n}\right),
\end{align*} for all $0< \varepsilon <1$, where $Z\sim N({\bf 0} ,I_n)$ and $M_{p,n}$ is a median for $\|Z\|_p$. 
Thus, for $p=5$ say, we get a contradiction for all $\varepsilon \in (n^{-1/6} ,1)$.

\subsection{Another small deviation inequality} In this subsection we prove a small deviation inequality similar to 
\eqref{eq:gsd} with the variance replaced by the $\mathbb E \|\nabla f\|_2^2$. The inequality is known to specialists and the method of 
proof goes at least back to \cite{BG} and \cite{Sam}. The inequality is weaker than \eqref{eq:gsd} but holds for a larger class of measures:
for all measures which satisfy a quadratic cost inequality \'a la Talagrand \cite{Tal-cost} (see also \cite{AS-cvx} for a recent 
development on the related subject). First we recall the necessary definitions. 
For any two Borel probability measures $\mu$ and $\nu$ on $\mathbb R^n$
the Wasserstein distance $W_2(\mu,\nu)$ is defined as 
\begin{align*}
W_2^2(\mu,\nu)=\inf_{\pi} \iint \|x-y\|_2^2 \, d\pi(x,y),
\end{align*} where the infimum is taken over all couplings (or matchings) $\pi$ of $\mu$ and $\nu$, i.e. $\pi$ 
has marginals $\mu$ and $\nu$ respectively. The Kullback-Leibler divergence (or relative entropy) of $\nu$ with respect to 
$\mu$ is defined by
\begin{align*}
D(\nu || \mu) = {\rm Ent}_\mu \left( \frac{d\nu}{d\mu} \right) =\int \log \frac{d\nu}{d\mu} \, d\nu,
\end{align*} if $\nu$ is absolutely continuous with respect to $\mu$ with Radon-Nikodym derivative $\frac{d\nu}{d\mu}$ and $\infty$ otherwise.

A Borel probability measure $\mu$ on $\mathbb R^n$ it is said to satisfy a {\it quadratic transportation cost inequality} with
constant $\rho >0$ if
\begin{align} \label{eq:transp-ineq}
W_2(\mu, \nu) \ls \sqrt{ \rho D( \nu || \mu) },
\end{align} for any Borel probability measure $\nu$ with $\nu \ll \mu$. It is known that measures with this property can be 
characterized in terms of infimum convolution inequalities with cost function $w(z)=\|z\|_2^2/(2\rho), \; z\in \mathbb R^n$ (e.g. see 
\cite[Corollary 6.4.]{Led}). Furthermore, Otto and Villani in \cite{OV} showed that measures which satisfy a log-Sobolev inequality 
they also satisfy a quadratic transportation cost inequality. The main inequality of this subsection, which is in the same spirit 
as Theorem \ref{thm:main-2}, reads as follows. 

\begin{theorem}
Let $\mu$ be any Borel probability measure on $\mathbb R^n$, which satisfies a quadratic transportation cost inequality \eqref{eq:transp-ineq}. 
Then, for any smooth, convex map $f:\mathbb R^n \to \mathbb R$ we have:
\begin{align*}
\log \mathbb E_\mu e^{-f} \ls -\mathbb E_\mu (f)  + \frac{\rho}{4} \mathbb E_\mu \|\nabla f\|_2^2.
\end{align*} In particular,
\begin{align*}
\mu \left(  \left\{ x: f(x)- \mathbb E_\mu(f) \ls -t \left(\mathbb E_\mu \|\nabla f\|_2^2 \right)^{1/2} \right\} \right) \ls e^{-t^2/\rho},
\end{align*} for all $t>0$.
\end{theorem}

\noindent {\it Proof.} Since $f$ is convex and smooth, for any $x,y\in \mathbb R^n$ we may write:
\begin{align*}
f(x)-f(y) \ls \langle \nabla f(x), x-y \rangle \ls \|\nabla f (x)\|_2 \|x-y\|_2.
\end{align*} Fix any probability measure $\nu$ with $\nu \ll \mu$ and let $\pi$ be any coupling of $\mu$ and $\nu$. Thus, integration
with respect to $\pi$ yields:
\begin{align*}
\mathbb E_\mu f - \mathbb E_\nu f \ls \int  \|\nabla f(x)\|_2 \cdot  \|x-y\|_2 \, d\pi(x,y) \ls \left(\mathbb E_\mu \|\nabla f\|_2^2\right)^{1/2}
\left( \int \|x-y\|_2^2 \, d\pi(x,y)\right)^{1/2}.
\end{align*} Since the left-hand side is fixed for any coupling $\pi$ of $\mu$ and $\nu$ we infer:
\begin{align} \label{eq:lap-bd}
\mathbb E_\mu f - \mathbb E_\nu f \ls 
\sqrt{\rho \mathbb E_\mu \|\nabla f\|_2^2 D(\nu || \mu) },
\end{align} where we have used the assumption on $\mu$. Now we employ Gibb's variational formula (for a proof see \cite[Corollary 4.14]{BLM}): 
For any $\mu$-measurable map $f$ one has
\begin{align}
\log \mathbb E_\mu e^f = \sup_{\nu \ll \mu} \left\{ \mathbb E_\nu f- D(\nu || \mu) \right\}.
\end{align} Applying the latter for $- f$ and taking into account \eqref{eq:lap-bd} we obtain:
\begin{align*}
\log \mathbb E_\mu e^{- f} \ls  - \mathbb E_\mu f +
\sup_{\nu \ll \mu} \left\{  \sqrt{\rho \mathbb E_\mu \|\nabla f\|_2^2 D(\nu ||\mu|)} -D (\nu || \mu) \right\} 
\ls -\mathbb E_\mu f +\frac{\rho}{4} \mathbb E_\mu \|\nabla f\|_2^2.
\end{align*} The result follows. \prend

\begin{remark}
In view of Theorem \ref{thm:sd-gauss} and Theorem \ref{thm:main-2} one might ask if in the above inequality we may replace 
the $(\mathbb E\|\nabla f\|_2^2)^{1/2}$ by the $\sqrt{{\rm Var}(f)}$ under the stronger assumption that the measure
$\mu$ satisfies a log-Sobolev inequality. Let us mention that this is not the case. For example the uniform measure $\mu_{D_n}$ 
on the Euclidean ball $D_n$ of volume 1, it is known to satisfy log-Sobolev inequality (e.g. see \cite{BobLed}) 
but the correct estimate for the small deviation of the norm $\|\cdot\|_{D_n}$ with respect to this measure is subexponential; 
see Remark \ref{rem:optimality-beta}.
\end{remark}

\medskip
\section{Probabilistic estimates on the Grassmannian}

In this Section we prove Theorem \ref{thm: conc-M}. We present two approaches to derive this probabilistic estimate. The first one
uses Gaussian tools and is based on the new small deviation inequality \eqref{eq:gsd}, therefore yields better tails estimate
in the one-sided small deviation regime, but restricts the range of $t$ and of dimension $k$. The second approach overcomes this obstacle
by working directly on the Grassmann space, but the tail estimate we obtain relies on the Lipschitz constant since
we employ the Gromov-Milman theorem for $SO(n)$. A small ball probability estimate for the mapping 
$F\mapsto w(P_FA)$ is also provided.

\subsection{From Gauss' space to Grassmannian}

First we provide a Gaussian proof of Theorem \ref{thm: conc-M}. Let us recall the following:

\begin{lemma} \label{lem:a1}
Let $A$ be symmetric convex body on $\mathbb R^n$. For any matrix $T=(t_{ij})_{i,j=1}^{k,n}, \, 1\ls k \ls n$ with rank $k$ we consider 
the map $T\mapsto w(TA)$. Then, we have the Lipschitz condition:
\begin{align} \label{eq:lip-w}
|w(TA)-w(SA)| \ls \frac{ R(A)}{\sqrt{k}} \|T-S\|_{\rm HS},
\end{align} for all $T=(t_{ij}),S=(s_{ij}) \in \mathbb R^{k\times n}$. Therefore, we have:
\begin{align} \label{eq:ld-w}
\mathbb P \left( w(GA) \gr \mathbb E [w(GA)] +t \right) \ls \exp(-ct^2 k /R(A)^2 ),
\end{align} and 
\begin{align} \label{eq:sd-w}
\mathbb P\left( w(GA) \ls \mathbb E[w(GA)] -t \right) \ls \exp \left( -ct^2 \max \left\{ \frac{1}{ {\rm Var}[h_A(Z)]} , \frac{k}{R(A)^2} \right\} \right), 
\end{align} for all $t>0$, where $G=(g_{ij})$ is $k\times n$ matrix with i.i.d. standard Gaussian entries and 
$Z\sim N({\bf 0},I_n)$.
\end{lemma}

\noindent {\it Proof.} The proof of \eqref{eq:lip-w} is standard and is left to the reader. Then, estimate \eqref{eq:ld-w}
immediately follows from the concentration on Gauss' space. For proving \eqref{eq:sd-w} we need the next fact:

\smallskip

\noindent {\it Claim.} If $G=(g_{ij})_{i,j=1}^{k,n}$ is a Gaussian matrix then, 
\begin{align} \label{eq:â-mw}
{\rm Var}[w(GA)] \ls \min \left\{\frac{R(A)^2}{k}, {\rm Var}[h_A(Z)] \right\},
\end{align} where $Z\sim N({\bf 0},I_n)$.

\smallskip

\noindent {\it Proof of Claim.} The bound in terms of the circumradius follows from the Lipschitz condition and the 
Poincar\'e inequality \eqref{eq:Poin}. For bounding in terms of the variance we use the Cauchy-Schwarz inequality.

Finally, the estimate \eqref{eq:sd-w} follows from the small deviation inequality \eqref{eq:gsd} applied for the convex function 
$T\mapsto w(TA)$. \prend

\medskip

One more ingredient is the polar decomposition of any matrix $T\in \mathbb R^{k\times n}$. If $T\in \mathbb R^{k\times n}$ 
we may write: $T= SQ$ where $S=(TT^\ast)^{1/2}$ and $Q$ is the orthogonal projection onto $F={\rm Im}T^\ast$.

The next Lemma follows if we take into account the above decomposition and the ideal property of the $\ell$-norm (see e.g. \cite{TJ}).

\begin{lemma} \label{lem:a2}
If $T\in \mathbb R^{k\times n}$ and $A$ is a symmetric convex body on $\mathbb R^n$, then we have:
\begin{align*}
\lambda_k( (TT^\ast)^{1/2} ) w(P_FA) \ls w(TA) \ls \lambda_1((TT^\ast)^{1/2}) w(P_FA),
\end{align*} where $\lambda_j( (TT^\ast)^{1/2}) \equiv s_j(T)$ is the $j$-th eigenvalue of $(TT^\ast)^{1/2}$ (or the $j$-th singular value of $T$) and $F={\rm Im}T^\ast$.
\end{lemma}

\noindent {\it Proof.} First note that if $S:\ell_2^k\to \ell_2^k$ is a linear map which satisfies $0<a \ls \|S\theta\|_2 \ls b$ for all $\theta\in S^{k-1}$, then 
\begin{align} \label{eq:Ap-2-side}
a \mathbb E[h_A(Y)] \ls \mathbb E [h_A(SY)] \ls b \mathbb E [h_A(Y)]
\end{align} where $Y\sim N({\bf 0},I_k)$. This follows by the ideal property of the $\ell$-norm, i.e. for any operator 
$u:\ell_2^n \to \ell_2^n$ and $v:\ell_2^n\to X$ we have $\ell(vu)\ls \ell(v) \|u\|_{\rm op}$.
Now in our setting, set $S=(TT^\ast)^{1/2}$. If $Y\sim N({\bf 0},I_k)$ and $c_k:= \mathbb E\|Y\|_2$, we may write:
\begin{align*}
w(TA) = c_k^{-1} \mathbb E [h_{TA}(Y)] = c_k^{-1} \mathbb E [h_{P_FA}(S^\ast Y)] \ls c_k^{-1} \|S\|_{2\to 2} \mathbb E [h_{P_FA}(Y)],
\end{align*} where we have used the right-hand side of \eqref{eq:Ap-2-side}. We work similarly for the lower estimate. \prend

\bigskip

In the Gaussian random setting the variables $\lambda_j((GG^\ast)^{1/2})$ and $w(P_FA)$ with $F= {\rm Im} G^\ast$ are independent each other. Namely, we have the following:

\begin{lemma} \label{lem:a3}
Let $A$ be a symmetric convex body on $\mathbb R^n$ and let $G=(g_{ij})_{i,j=1}^{k,n}$ where $g_{ij}$ are i.i.d. standard normals. Then, $F={\rm Im}G^\ast$ is uniformly distributed over $G_{n,k}$ and the random variables $\lambda_j( (GG^\ast)^{1/2})$ and
$w(P_FA)$ are independent.
\end{lemma}

The proof of this fact follows the same lines as in \cite[Proposition 4.1]{PPZ}. Now we are ready to prove the following:

\begin{theorem}
Let $A$ be a symmetric convex body in $\mathbb R^n$. Fix $1\ls k\ls n-1$. Then, we have:
\begin{align*}
\nu_{n,k} \left( F\in G_{n,k} : w(P_FA) > (1+t) w(A) \right) \ls C\exp(-ct^2 k k_\ast(A) ), 
\end{align*} for all $t> c_1\sqrt{k/n}$. Furthermore,
\begin{align*}
\nu_{n,k} \left( F\in G_{n,k} : w(P_FA) \ls (1-t) w(A) \right) \ls C \exp\left( -ct^2 \max \left\{ k k_\ast(A), \frac{1}{\beta_\ast(A)} \right \} \right),
\end{align*} for all $ c_2\sqrt{k/n} <t<1$, provided that $k\lesssim n$.
\end{theorem}

\noindent {\it Proof.} Note that for any Gaussian matrix $G=(g_{ij})_{i,j=1}^{k,n}$ we have: 
\begin{align*}
\mathbb E [w(GA)] = \mathbb E [h_A(Z)] = \mathbb E\|Z\|_2 \cdot w(A),
\end{align*}
where $Z\sim N({\bf 0},I_n)$. Fix $t>0$. From Lemma \ref{lem:a2} and Lemma \ref{lem:a3} we may write:
\begin{align*}
 \nu_{n,k}( w(P_FA) >(1+t)w(A)) \mathbb P\left(\lambda_k((GG^\ast)^{1/2} )> (1-\delta) \mathbb E\|Z\|_2 \right) \\
 \ls \mathbb P( w(GA) > (1+t/2)\mathbb E w(GA)),
\end{align*} for $0<\delta<\frac{t}{2(1+t)}$. Recall the following well known:

\smallskip

\noindent {\it Fact.} Let $\delta\in (0,1)$. Then the random Gaussian matrix $G=(g_{i,j})_{i,j=1}^{k,n}$ with $k\ls c\delta^2n$ satisfies:
\begin{align*}
(1-\delta) \mathbb E\|Z\|_2 < \lambda_k ((GG^\ast)^{1/2}) \ls \lambda_1( (GG^\ast)^{1/2}) < (1+\delta) \mathbb E\|Z\|_2
\end{align*} with probability greater than $1-e^{-c\delta^2 n}$, where $Z\sim N({\bf 0}, I_n)$. 

\smallskip

Now Lemma \ref{lem:a1} yields:
\begin{align*}
(1-e^{-c\delta^2 n}) \nu_{n,k} \left( w(P_FA) > (1+t) w(A) \right) \ls C\exp \left(-ct^2 k k_\ast(A) \right).
\end{align*} The choice $\delta\simeq \sqrt{k/n}$ yields the upper estimate. We work similarly for the lower estimate. \prend

\subsubsection{A small ball estimate} Next, we prove the following:

\begin{theorem} [small ball for the mean width of projections]
Let $A$ be a symmetric convex body in $\mathbb R^n$ and let $1\ls k\ls n-1$. Then, we have:
\begin{align*}
\nu_{n,k} \left( \left\{ F\in G_{n,k} : w(P_FA) \ls c\varepsilon w(A) \right\} \right) \ls (C \varepsilon)^{c \max \{ k k_\ast(A), \frac{1}{\beta_\ast(A)} \} },
\end{align*} for all $\varepsilon \in (0,1/2)$.
\end{theorem}

\noindent {\it Proof.} First note that the function $T\mapsto w(TA)$ is indeed a norm on $\mathbb R^{k\times n}$. 
If $\mathcal C_A =\left\{ T\in \mathbb R^{k\times n}: w(TA)\ls 1 \right\}$ is its unit ball, then estimate \eqref{eq:â-mw} 
shows that
$\beta( \mathcal C_A) \ls \min\left\{ \frac{1}{k k_\ast(A)}, \beta_\ast(A) \right \}$. It follows by \cite[Proposition 3.2]{PVsd} that:
\begin{align*}
\mathbb P \left( w(GA) \ls c\varepsilon \mathbb E[w(GA)] \right) \ls \frac{1}{2} \varepsilon^{c/\beta(\mathcal C_A)},
\end{align*} for all $\varepsilon \in (0,1/2)$. Now we use Lemma \ref{lem:a2} and Lemma \ref{lem:a3} to get:
\begin{align*}
c_1 \nu_{n,k} \left( w(P_FA) \ls c' \varepsilon w(A) \right) &\ls \nu_{n,k} \left( w(P_FA) \ls c' \varepsilon w(A) \right) 
\mathbb P \left( \lambda_1((GG^\ast)^{1/2}) \ls C_1 \mathbb E\|Y\|_2 \right) \\
&\ls \mathbb P \left( w(GA) \ls c\varepsilon \mathbb E[w(GA)] \right).
\end{align*} The result readily follows. \prend

\bigskip

\subsection{Breaking the barrier} Note that the argument of the preceding proof does not allow to consider 
neither subspaces of small codimension nor very small values of $t>0$. In order to adjust this technicality
we have to work beyond the Gaussian setting. First we state Lipschitz
estimates on the Grassmannian for the map $F\mapsto w(P_FA)$ with respect to the normalized metrics:
\begin{align*}
\sigma_\infty(E,F) = \|P_E-P_F\|_{\rm op}, \quad \sigma_2(E, F)= \frac{\|P_E - P_F\|_{\rm HS}}{\sqrt{k}}.
\end{align*}

\begin{lemma} \label{lem:lip-w-g}
  Let $A$ be a (symmetric) convex body in $\mathbb R^n$ and
  let $1\ls k\ls n-1$. Then, for any $E,F\in G_{n,k}$ we have:
  \begin{align}
    |w(P_EA)-w(P_FA)|\ls c \sqrt{n/k}w(A)\sigma_\infty(E,F).
  \end{align} Furthermore, we have:
  \begin{align}
    |w(P_EA)-w(P_FA)|\ls c'R(A)\sigma_2(E,F),
  \end{align} where $c,c'>0$ are absolute constants. Equivalently, in the dual setting, we may write:
  \begin{align*}
   |M(A\cap E)-M(A\cap F)| \lesssim \min \left\{ \sqrt{n/k}M(A)\sigma_\infty(E,F), b(A)\sigma_2(E,F)\right\}.
  \end{align*}

\end{lemma}

\noindent {\it Proof.} Using the formula 
\begin{align} \label{eq:pc}
c_k w(P_FA)=c_k\int_{S_F} h_A(\phi)\, d\sigma_F(\phi)=\int_F h_A(x)\, d\gamma_F(x)=\int_{\mathbb R^n} h_A(P_Fx)\, d\gamma_n(x),
\end{align} where $c_k= \mathbb E\|Y\|_2=\sqrt{2}\Gamma(\frac{k+1}{2})/\Gamma(\frac{k}{2}) \simeq \sqrt{k}$ 
with $Y\sim N({\bf 0},I_k)$, we may write:
\begin{align*}
  c_k|w(P_EA)-w(P_FA)| &\ls \int_{\mathbb R^n} |h_A(P_Ex)-h_A(P_Fx)|\, d\gamma_n(x) \\
  &\ls \int_{\mathbb R^n} h_A(P_Ex-P_Fx)\, d\gamma_n(x). 
\end{align*} Now we proceed as follows. In order to prove the first estimate recall the ideal property of the
$\ell$-norm (e.g. see \cite{TJ}). Applying this for $u=P_E-P_F$, $v=id$ and $X=(\mathbb R^n, \|\cdot\|_{A^\circ})$ we obtain:
\begin{align*}
 \int_{\mathbb R^n} h_A((P_E-P_F)x)\, d\gamma_n(x)\ls \int_{\mathbb R^n} h_A(x) \, d\gamma_n(x) \cdot \|P_E-P_F\|_{\rm op}.
\end{align*}
Therefore we get:
\begin{align*}
 |w(P_EA) - w(P_FA)|\ls \frac{c_n}{c_k}w(A) \|P_E-P_F\|_{\rm op}.
\end{align*}

\noindent For the second assertion we proceed as follows:

\begin{align*}
  \int_{\mathbb R^n} h_A(P_Ex-P_Fx)\, d\gamma_n(x) &\ls R(A) \int_{\mathbb R^n} \|(P_E-P_F)x\|_2\, d\gamma_n(x) \\
  &\ls R(A) \|P_E-P_F\|_{\rm HS},
\end{align*} by the Cauchy-Schwarz inequality and the isotropicity of $\gamma_n$.  \prend

\medskip

Let us recall the 
concentration on $SO(n)$. The next result is due to Gromov and V. Milman from \cite{GM}.

\begin{theorem} \label{thm:GroMil}  
Let $f : SO(n)\to \mathbb R$ be $L$-Lipschitz map with respect to the Hilbert-Schmidt
 norm, i.e. $|f(U)-f(V)|\ls L \|U-V\|_{\rm HS}$ for any $U,V\in SO(n)$. Then we have:
 \begin{align*}
  \mu_n \left( \{U\in SO(n) : |f(U)-\mathbb Ef|\gr t\}\right)\ls C_1\exp(-c_1 nt^2/L^2),
 \end{align*} for all $t>0$.
\end{theorem}

First note that the distribution of the map $F\mapsto w(P_FA)$ can be carried out over the orthogonal 
group $O(n)$. This follows directly
 from the Haar-invariance of $\nu_{n,k}$ under $O(n)$. Indeed,
\begin{align*}
  \nu_{n,k}(\{F\in G_{n,k} : w(P_FA)\in B\}) &= \mu_n(\{ U\in O(n) : w(P_{UE}A)\in B\})
\end{align*} for some (any) fixed $E\in G_{n,k}$. That is
 $F\mapsto w(P_FA)$ can be equivalently viewed as function of $U\in O(n)$. So, in order to apply a concentration result for $F\mapsto w(P_FA)$ on $G_{n,k}$ it suffices to suitably apply the Gromov-Milman theorem for $U\mapsto w(P_{UE}A)$ for some fixed 
$E\in G_{n,k}$. Toward this end we shall need the Lipschitz constant of $U \mapsto \psi(U)= w(P_{UE}A)$ with respect to the 
Hilbert-Schmidt norm. Then Lemma \ref{lem:lip-w-g} implies:
\begin{align*}
  |\psi(U_1)- \psi(U_2)| \lesssim \frac{R(A)}{ \sqrt{k} } \| P_{U_1E}-P_{U_2E} \|_{\rm HS}.
\end{align*} On the other hand we have\footnote{Set $U=U_1^\ast U_2$. Using the (right) invariance under the orthogonal group and the 
contractive property of the Hilbert-Schmidt norm, we may write: 
\begin{align*}
\| P_{U_1E}-P_{U_2E} \|_{\rm HS} = \|P_E - UP_EU^\ast \|_{\rm HS} \ls \| P_E(I-U^\ast)\|_{\rm HS} + \|(I-U) P_E U^\ast \|_{\rm HS} 
&\ls \| I-U^\ast \|_{\rm HS} + \| (I-U) P_E \|_{\rm HS} \\
&\ls 2 \| I- U \|_{\rm HS}= 2\|U_1-U_2\|_{\rm HS}.\end{align*} } 
$ \| P_{U_1E}-P_{U_2E} \|_{\rm HS} \ls 2 \|U_1-U_2\|_{\rm HS}$, which proves the Lipschitz
condition for $\psi$.

Now we explain how we derive the concentration estimate for $\psi$ on the full orthogonal group by following an argument 
from \cite{Me}. Let $U_1$ be Haar distributed
on $SO(n)$ and let $M_\pi$ be the permutation matrix\footnote{For any permutation $\pi$ of the set $\{1,2,\ldots,n\}$ the 
permutation matrix $M_\pi$ acts as: $M_\pi(e_j)=e_{\pi(j)}$ for $j=1,\ldots,n$.} which corresponds to the permutation 
$\pi$ such that $\pi(j)=j$ for $j\ls n-2$, $\pi (n-1)=n$ and
$\pi (n)=n-1$. Then, $U_2:=U_1M_\pi$ is the same as $U_1$ (as matrices) but 
 the last two rows which are switched. Therefore, $U_2$ is Haar distributed on $SO^-(n)$. Define the orthogonal map $U$ being $U_1$ with probability
 1/2 and $U_2$ with probability 1/2. It follows that $U$ is Haar distributed over $O(n)$. Moreover, note that 
 since $E$ can be considered as $E=[e_i : 1\ls i\ls k]$ we have $M_\pi(E)=E$ as long as $k\ls n-2$. In that case:
 \begin{align*} 
 \mathbb E_U\psi(U) = \mathbb E_{U_1} \psi(U_1) =\mathbb E_{U_2} \psi(U_2) =\int_{G_{n,k}} w(P_FA) \, d\nu_{n,k}(F)=w(A).
 \end{align*} Conditioning on whether $U=U_1$ or $U=U_2$ and taking into account Theorem \ref{thm:GroMil} we get:
 \begin{align*}
 \mu_n \left( \left\{ U\in O(n) : \left| \psi(U) - \int_{O(n)} \psi \, d\mu_n \right| >t \right\} \right) \ls c_1\exp(-c_2 n t^2 k/R(A)^2),
 \end{align*} for all $t>0$. The preceding discussion leads us to the next:

\begin{theorem}[concentration for the mean width of projections] \label{thm:conc-mw}
  Let $A$ be a symmetric convex body on $\mathbb R^n$ and let $1\ls k\ls n-1$.
  Then, one has the following concentration inequality:
  \begin{align*}
    \nu_{n,k}\left( \big\{ F\in G_{n,k} : |w(P_FA)-w(A)|\gr t w(A) \big \} \right) \ls c_1\exp \left( -c_2 t^2 k k_\ast(A) \right),
  \end{align*} for all $t>0$, where $c_1,c_2>0$ are absolute constants.
\end{theorem}

\noindent {\it Proof.} The case $1\ls k\ls n-2$ follows from the previous argument. We treat the case $k=n-1$ separately. Note that any $F\in G_{n,n-1}$ can be identified with $e^\perp$ for 
$e \in S^{n-1}$. Thus, we may write:
\begin{align*}
\int_{G_{n,n-1}} w(P_FA) \, d\nu_{n,n-1}(F) = \int_{S^{n-1}} w(P_{e^\perp}A) \, d\sigma(e) = w(A),
\end{align*} and as in \eqref{eq:pc} we have:
\begin{align*}
w(P_{e^\perp} A)= \frac{c_n}{c_{n-1}} \int_{S^{n-1}} h_A(P_{e^\perp} \theta) \, d\sigma (\theta) = \frac{c_n}{c_{n-1}} \int_{S^{n-1}} h_A(\theta-\langle \theta ,e \rangle e) \, d\sigma(\theta).
\end{align*} Since $c_n/c_{n-1}\simeq 1$, it suffices to work with the map $f:S^{n-1}\to \mathbb R$ defined by $f(e):= \int_{S^{n-1}} h_A(\theta-\langle \theta, e \rangle e) \, d\sigma(\theta)$. We have the following:

\smallskip

\noindent {\it Claim.} The mapping $f:S^{n-1}\to \mathbb R$ is Lipschitz with:
\begin{align*}
{\rm Lip}(f) \ls R(A) \sqrt{\frac{2}{n}}.
\end{align*}

\smallskip

\noindent {\it Proof of Claim.} For any $u,v\in S^{n-1}$, we have:
\begin{align*}
|f(u)-f(v)| &\ls \int_{S^{n-1}} h_A\left( \langle u,\theta\rangle u - \langle v,\theta \rangle v \right) \, d\sigma(\theta) \\
&\ls R(A) \left( \int_{S^{n-1}}\left \| \langle u,\theta \rangle u- \langle v,\theta \rangle v \right\|_2^2 \, d\sigma(\theta) \right)^{1/2}.
\end{align*} Note that:
\begin{align*}
\int_{S^{n-1}}\left \| \langle u,\theta \rangle u- \langle v,\theta \rangle v \right\|_2^2 \, d\sigma(\theta) &= 2\int_{S^{n-1}} \theta_1^2\, d\sigma(\theta) - 2\langle u,v\rangle \int_{S^{N-1}} \langle u,\theta \rangle \langle v,\theta \rangle \, d\sigma(\theta) \\
&= \frac{2}{n} (1-\langle u,v\rangle^2) \ls \frac{2}{n} \|u-v\|_2^2.
\end{align*}
The result now follows from the concentration on the sphere $S^{n-1}$. \prend

\bigskip

\noindent {\bf Acknowledgements.} The authors are grateful to Peter Pivovarov for useful remarks. They would also like to thank the
anonymous referee whose helpful comments improved the presentation of the paper.

\bigskip
\bibliography{var-euclid-ref}

\begin{thebibliography}{AAGM15}

\bibitem[AAGM15]{AGM}
S.~Artstein-Avidan, A.~Giannopoulos, and V.~D. Milman.
\newblock {\em Asymptotic geometric analysis. {P}art {I}}, volume 202 of {\em
  Mathematical Surveys and Monographs}.
\newblock American Mathematical Society, Providence, RI, 2015.

\bibitem[AK16]{AK}
F.~Albiac and N.~Kalton.
\newblock {\em Topics in {B}anach space theory}, volume 233 of {\em Graduate
  Texts in Mathematics}.
\newblock Springer, [Cham], second edition, 2016.
\newblock With a foreword by Gilles Godefory.

\bibitem[AS17]{AS-cvx}
R.~Adamczak and C.~Strzeleki.
\newblock On the convex {P}oincar\'{e} inequality and weak transportation
  transportation inequalities.
\newblock \url{https://arxiv.org/abs/1703.01765}, (2017).
\newblock preprint.

\bibitem[BG99]{BG}
S.~G. Bobkov and F.~G\"otze.
\newblock Exponential integrability and transportation cost related to
  logarithmic {S}obolev inequalities.
\newblock {\em J. Funct. Anal.}, 163(1):1--28, 1999.

\bibitem[BL00]{BobLed}
S.~G. Bobkov and M.~Ledoux.
\newblock From {B}runn-{M}inkowski to {B}rascamp-{L}ieb and to logarithmic
  {S}obolev inequalities.
\newblock {\em Geom. Funct. Anal.}, 10(5):1028--1052, 2000.

\bibitem[BLM13]{BLM}
S.~Boucheron, G.~Lugosi, and P.~Massart.
\newblock {\em Concentration inequalities}.
\newblock Oxford University Press, Oxford, 2013.
\newblock A nonasymptotic theory of independence, With a foreword by Michel
  Ledoux.

\bibitem[Bob99]{Bob}
S.~G. Bobkov.
\newblock Isoperimetric and analytic inequalities for log-concave probability
  measures.
\newblock {\em Ann. Probab.}, 27(4):1903--1921, 1999.

\bibitem[Bog98]{Bog}
V.~I. Bogachev.
\newblock {\em Gaussian measures}, volume~62 of {\em Mathematical Surveys and
  Monographs}.
\newblock American Mathematical Society, Providence, RI, 1998.

\bibitem[Bor75]{Bor}
C.~Borell.
\newblock Convex set functions in {$d$}-space.
\newblock {\em Period. Math. Hungar.}, 6(2):111--136, 1975.

\bibitem[CEFM04]{CFM}
D.~Cordero-Erausquin, M.~Fradelizi, and B.~Maurey.
\newblock The ({B}) conjecture for the {G}aussian measure of dilates of
  symmetric convex sets and related proble.
\newblock {\em J. Funct. Anal.}, 214(2):410--427, 2004.

\bibitem[CEL12]{CorLed}
D.~Cordero-Erausquin and M.~Ledoux.
\newblock Hypercontractive measures, {T}alagrand's inequality, and influences.
\newblock In {\em Geometric aspects of functional analysis}, volume 2050 of
  {\em Lecture Notes in Math.}, pages 169--189. Springer, Heidelberg, 2012.

\bibitem[Cha14]{Cha}
S.~Chatterjee.
\newblock {\em Superconcentration and related topics}.
\newblock Springer Monographs in Mathematics. Springer, Cham, 2014.

\bibitem[DR50]{DR}
A.~Dvoretzky and C.~A. Rogers.
\newblock Absolute and unconditional convergence in normed linear spaces.
\newblock {\em Proc. Nat. Acad. Sci. U. S. A.}, 36:192--197, 1950.

\bibitem[Dvo61]{Dvo}
A.~Dvoretzky.
\newblock Some results on convex bodies and {B}anach spaces.
\newblock In {\em Proc. Internat. Sympos. Linear Spaces (Jerusalem, 1960)},
  pages 123--160. Jerusalem Academic Press, Jerusalem; Pergamon, Oxford, 1961.

\bibitem[Ehr83]{Ehr}
A.~Ehrhard.
\newblock Sym\'etrisation dans l'espace de {G}auss.
\newblock {\em Math. Scand.}, 53(2):281--301, 1983.

\bibitem[ENT16]{ENT}
A.~Eskenazis, P.~Nayar, and T.~Tkocz.
\newblock Gaussian mixtures: entropy and geometric inequalities.
\newblock \url{https://arxiv.org/abs/1611.04921}, (2016).
\newblock preprint.

\bibitem[Fig76]{F}
T.~Figiel.
\newblock A short proof of {D}voretzky's theorem on almost spherical sections
  of convex bodies.
\newblock {\em Compositio Math.}, 33(3):297--301, 1976.

\bibitem[FLM77]{FLM}
T.~Figiel, J.~Lindenstrauss, and V.~D. Milman.
\newblock The dimension of almost spherical sections of convex bodies.
\newblock {\em Acta Math.}, 139(1-2):53--94, 1977.

\bibitem[GM83]{GM}
M.~Gromov and V.~D. Milman.
\newblock A topological application of the isoperimetric inequality.
\newblock {\em Amer. J. Math.}, 105(4):843--854, 1983.

\bibitem[Gor85]{Go}
Y.~Gordon.
\newblock Some inequalities for {G}aussian processes and applications.
\newblock {\em Israel J. Math.}, 50(4):265--289, 1985.

\bibitem[HW17]{HW}
H.~Huang and F.~Wei.
\newblock Upper bound for the {D}voretzky dimension in {M}ilman-{S}chechtman
  theorem.
\newblock In {\em Geometric aspects of functional analysis}, volume 2169 of
  {\em Lecture Notes in Math.}, pages 181--186. Springer, Cham, 2017.

\bibitem[JL84]{JL}
W.~B. Johnson and J.~Lindenstrauss.
\newblock Extensions of {L}ipschitz mappings into a {H}ilbert space.
\newblock In {\em Conference in modern analysis and probability ({N}ew {H}aven,
  {C}onn., 1982)}, volume~26 of {\em Contemp. Math.}, pages 189--206. Amer.
  Math. Soc., Providence, RI, 1984.

\bibitem[Joh48]{Jo}
F.~John.
\newblock Extremum problems with inequalities as subsidiary conditions.
\newblock In {\em Studies and {E}ssays {P}resented to {R}. {C}ourant on his
  60th {B}irthday, {J}anuary 8, 1948}, pages 187--204. Interscience Publishers,
  Inc., New York, N. Y., 1948.

\bibitem[Kla04]{Kl1}
B.~Klartag.
\newblock A geometric inequality and a low {$M$}-estimate.
\newblock {\em Proc. Amer. Math. Soc.}, 132(9):2619--2628 (electronic), 2004.

\bibitem[Kla07]{Kl-sub}
B.~Klartag.
\newblock Uniform almost sub-{G}aussian estimates for linear functionals on
  convex sets.
\newblock {\em Algebra i Analiz}, 19(1):109--148, 2007.

\bibitem[KV07]{KV}
B.~Klartag and R.~Vershynin.
\newblock Small ball probability and {D}voretzky's theorem.
\newblock {\em Israel J. Math.}, 157:193--207, 2007.

\bibitem[Led01]{Led}
M.~Ledoux.
\newblock {\em The concentration of measure phenomenon}, volume~89 of {\em
  Mathematical Surveys and Monographs}.
\newblock American Mathematical Society, Providence, RI, 2001.

\bibitem[Lew78]{Lew}
D.~R. Lewis.
\newblock Finite dimensional subspaces of {$L_{p}$}.
\newblock {\em Studia Math.}, 63(2):207--212, 1978.

\bibitem[LO05]{LO}
R.~Lata{\l}a and K.~Oleszkiewicz.
\newblock Small ball probability estimates in terms of widths.
\newblock {\em Studia Math.}, 169(3):305--314, 2005.

\bibitem[LT11]{LT}
M.~Ledoux and M.~Talagrand.
\newblock {\em Probability in {B}anach spaces}.
\newblock Classics in Mathematics. Springer-Verlag, Berlin, 2011.
\newblock Isoperimetry and processes, Reprint of the 1991 edition.

\bibitem[Mec14]{Me}
E.~Meckes.
\newblock Concentration of measure on the compact classical matrix groups.
\newblock \url{http://www.case.edu/artsci/math/esmeckes/Haar_notes.pdf},
  (2014).
\newblock Lecture Notes.

\bibitem[Mil71]{Mil1}
V.~D. Milman.
\newblock A new proof of {A}. {D}voretzky's theorem on cross-sections of convex
  bodies.
\newblock {\em Funkcional. Anal. i Prilo\v zen.}, 5(4):28--37, 1971.

\bibitem[MP89]{MiPa}
V.~D. Milman and A.~Pajor.
\newblock Isotropic position and inertia ellipsoids and zonoids of the unit
  ball of a normed {$n$}-dimensional space.
\newblock In {\em Geometric aspects of functional analysis (1987--88)}, volume
  1376 of {\em Lecture Notes in Math.}, pages 64--104. Springer, Berlin, 1989.

\bibitem[MS86]{MS}
V.~D. Milman and G.~Schechtman.
\newblock {\em Asymptotic theory of finite-dimensional normed spaces}, volume
  1200 of {\em Lecture Notes in Mathematics}.
\newblock Springer-Verlag, Berlin, 1986.
\newblock With an appendix by M. Gromov.

\bibitem[MS97]{MS-global}
V.~D. Milman and G.~Schechtman.
\newblock Global versus local asymptotic theories of finite-dimensional normed
  spaces.
\newblock {\em Duke Math. J.}, 90(1):73--93, 1997.

\bibitem[OV00]{OV}
F.~Otto and C.~Villani.
\newblock Generalization of an inequality by {T}alagrand and links with the
  logarithmic {S}obolev inequality.
\newblock {\em J. Funct. Anal.}, 173(2):361--400, 2000.

\bibitem[Pis89]{Pis_book}
G.~Pisier.
\newblock {\em The volume of convex bodies and {B}anach space geometry},
  volume~94 of {\em Cambridge Tracts in Mathematics}.
\newblock Cambridge University Press, Cambridge, 1989.

\bibitem[PPZ14]{PPZ}
G.~Paouris, P.~Pivovarov, and J.~Zinn.
\newblock A central limit theorem for projections of the cube.
\newblock {\em Probab. Theory Related Fields}, 159(3-4):701--719, 2014.

\bibitem[PV15]{PV-dvo-L_p}
G.~Paouris and P.~Valettas.
\newblock On {D}voretzky's theorem for subspaces of ${L}_p$.
\newblock \url{http://arxiv.org/abs/1510.07289}, (2015).
\newblock preprint.

\bibitem[PV16]{PVsd}
G.~Paouris and P.~Valettas.
\newblock A {G}aussian small deviation inequality for convex functions.
\newblock Annals of Probability (to appear), 2016.

\bibitem[PVZ17]{PVZ}
G.~Paouris, P.~Valettas, and J.~Zinn.
\newblock Random version of {D}voretzky's theorem in {$\ell_p^n$}.
\newblock {\em Stochastic Process. Appl.}, 127(10):3187--3227, 2017.

\bibitem[Sam03]{Sam}
P.~M. Samson.
\newblock Concentration inequalities for convex functions on product spaces.
\newblock In {\em Stochastic inequalities and applications}, volume~56 of {\em
  Progr. Probab.}, pages 33--52. Birkh\"auser, Basel, 2003.

\bibitem[Sch89]{Sch1}
G.~Schechtman.
\newblock A remark concerning the dependence on {$\epsilon$} in {D}voretzky's
  theorem.
\newblock In {\em Geometric aspects of functional analysis (1987--88)}, volume
  1376 of {\em Lecture Notes in Math.}, pages 274--277. Springer, Berlin, 1989.

\bibitem[Sch07]{Sch-cube}
G.~Schechtman.
\newblock The random version of {D}voretzky's theorem in {$\ell^n_\infty$}.
\newblock In {\em Geometric aspects of functional analysis}, volume 1910 of
  {\em Lecture Notes in Math.}, pages 265--270. Springer, Berlin, 2007.

\bibitem[SZ90]{SchZ}
G.~Schechtman and J.~Zinn.
\newblock On the volume of the intersection of two {$L^n_p$} balls.
\newblock {\em Proc. Amer. Math. Soc.}, 110(1):217--224, 1990.

\bibitem[Sza74]{Sz}
A.~Szankowski.
\newblock On {D}voretzky's theorem on almost spherical sections of convex
  bodies.
\newblock {\em Israel J. Math.}, 17:325--338, 1974.

\bibitem[Tal94]{Tal-russo}
M.~Talagrand.
\newblock On {R}usso's approximate zero-one law.
\newblock {\em Ann. Probab.}, 22(3):1576--1587, 1994.

\bibitem[Tal96]{Tal-cost}
M.~Talagrand.
\newblock Transportation cost for {G}aussian and other product measures.
\newblock {\em Geom. Funct. Anal.}, 6(3):587--600, 1996.

\bibitem[Tik14]{Tik-cube}
K.~E. Tikhomirov.
\newblock The randomized {D}voretzky's theorem in {$l_\infty^n$} and the
  {$\chi$}-distribution.
\newblock In {\em Geometric aspects of functional analysis}, volume 2116 of
  {\em Lecture Notes in Math.}, pages 455--463. Springer, Cham, 2014.

\bibitem[Tik17]{Tik-1-unc}
K.~E. Tikhomirov.
\newblock Superconcentration, and randomized {D}voretzky's theorem for spaces
  with 1-unconditional bases.
\newblock {\em Journal of Functional Analysis}, (to appear), 2017.
\newblock \url{https://doi.org/10.1016/j.jfa.2017.08.021}.

\bibitem[TJ89]{TJ}
N.~Tomczak-Jaegermann.
\newblock {\em Banach-{M}azur distances and finite-dimensional operator
  ideals}, volume~38 of {\em Pitman Monographs and Surveys in Pure and Applied
  Mathematics}.
\newblock Longman Scientific \& Technical, Harlow; copublished in the United
  States with John Wiley \& Sons, Inc., New York, 1989.

\bibitem[Val17]{V-tight}
P.~Valettas.
\newblock On the tightness of {G}aussian concentration for convex functions.
\newblock \url{https://arxiv.org/abs/1706.09446}, (2017).
\newblock preprint.

\end{thebibliography}
\bibliographystyle{alpha}

\end{document}